\documentclass{amsart}
\usepackage{amsmath}
\usepackage{amssymb}
\usepackage{amscd}
\usepackage{graphicx}
\usepackage{latexsym}
\usepackage{psfrag}
\usepackage{framed, color}

\usepackage{mathtools}

\usepackage{doc}

\newtheorem{thm}{Theorem}[section]
\newtheorem{lem}[thm]{Lemma}

\theoremstyle{definition}
\newtheorem{defn}[thm]{Definition}

\theoremstyle{remark}
\newtheorem{rk}[thm]{Remark}

\newtheorem{problem}[thm]{Problem}

\makeatletter
 
 \@addtoreset{equation}{section}
\makeatother

\newcommand{\Int}{\mathop{\mathrm{Int}}\nolimits}

\newcommand{\R}{\mathbb{R}}
\newcommand{\Z}{\mathbb{Z}}

\newcommand{\bang}{\textbf{\Large !\,}}

\def\spmapright#1{\smash{%
 \mathop{\hbox to 1.3cm{\rightarrowfill}}
  \limits^{#1}}}
\def\spmapleft#1{\smash{%
 \mathop{\hbox to 1.3cm{\leftarrowfill}}
  \limits^{#1}}}

\begin{document}

\title[Simplifying generic maps to the $2$--sphere and to the plane]{Simplifying 
generic smooth maps to the $2$--sphere and to the plane}

\author[O.~Saeki]{Osamu Saeki}
\address{Institute of Mathematics for Industry,
Kyushu University, Motooka 744, Nishi-ku, Fukuoka 819-0395,
Japan}

\date{\today}

\begin{abstract} 
We study how to construct explicit deformations of generic smooth
maps from closed $n$--dimensional
manifolds $M$ with $n \geq 2$ to the $2$--sphere $S^2$ and show that every smooth
map $M \to S^2$ is homotopic to a $C^\infty$ stable map with at most one cusp point
and with only folds of the middle absolute index. Furthermore, if $n$ is even, such a
$C^\infty$ stable map can be so constructed that the restriction to the singular point
set is a topological embedding. As a corollary, we show that for $n \geq 2$ even,
there always exists a $C^\infty$ stable map $M \to \R^2$ with at most one
cusp point such that the restriction to the singular point set
is a topological embedding. As another corollary, we give
a new proof to the existence of an open book structure
on odd dimensional manifolds which extends a given
one on the boundary, originally due to Quinn.
Finally, using the open book structure thus constructed,
we show that $k$--connected
$n$--dimensional manifolds always admit a fold map
into $\R^2$ without folds of absolute indices $i$ 
with $1 \leq i \leq k$,
for $n \geq 7$ odd and $1 \leq k \leq (n-5)/2$.
\end{abstract}

\maketitle

\section{Introduction}\label{section1}

For an arbitrary smooth closed manifold $M$ of dimension $n \geq 1$,
we have a Morse function $M \to \R$.
It is known that by modifying such a function by homotopy
if necessary, we can arrange
so that the critical values are in the 
order of indices, i.e.\ 
for every pair $p, q$ of critical points of indices $i(p), i(q)$, 
respectively, with
$i(p) < i(q)$, we have $f(p) < f(q)$, where $f : M \to \R$
denotes the modified Morse function.

Furthermore, Smale \cite{Sm2} showed that
if an $n$--dimensional closed manifold
$M$ with $n \geq 6$ is $k$--connected with $k \geq 1$, then
there exists a Morse function $f: M \to \R$
without critical points of indices $i$ with $1 \leq i \leq k$
and $n-k \leq i \leq n-1$. This was a very important
result in differentiable topology, which showed
a strong geometric consequence of an algebraic topological
property for manifolds, and which led to the solution
to the Poincar\'e Conjecture in high dimensions.

So, the following problem naturally arises.
For the definitions and properties of generic maps
or $C^\infty$ stable maps into manifolds of dimension $2$, 
see \S\ref{section2}.

\begin{problem}\label{prob1}
How about generic maps or $C^\infty$ stable maps
to manifolds of dimensions $\geq 2$?
The singular values can be nicely arranged with respect
to the indices?
If the manifold is highly connected, then the
singular points of certain indices can be eliminated?
\end{problem}

We note that such a problem has been considered in \cite{GK2},
for example. Recall that for a $C^\infty$ stable map $f$
of an $n$--dimensional manifold $M$, $n \geq 2$, into a
manifold of dimension $2$, we have only folds and cusps
as its singularities, and each of them has an absolute
index (for details, see \S\ref{section2}), which
is a concept similar to that for a nondegenerate
critical point of a smooth function. Furthermore,
the singular point set is a closed $1$--dimensional
submanifold of $M$.

In the present paper, we consider the above
problem for generic smooth
maps to $S^2$ and to $\R^2$.
Our first result is Theorem~\ref{thm1} (see also Remark~\ref{remarkn2}):
for a closed connected $n$--dimensional manifold $M$ with $n \geq 2$,
there always exists a $C^\infty$ stable
map $f : M \to S^2$ such that
for $n$ odd, $f$ has only folds of 
absolute index $(n-1)/2$ and has no cusps
and that for $n$ even, 
$f$ has only folds of absolute
index $(n-2)/2$ and at most one cusp, whose absolute index is equal to
$(n-2)/2$. Furthermore, when $n$ is even, we can choose $f$
such that the restriction of $f$
to the singular point set $S(f)$ is a topological embedding.
In the case of $n$ even, using the latter result, we will show,
as Theorem~\ref{thm2} (see also Remark~\ref{remarkn22}), that
there always exists a $C^\infty$ stable 
map $f : M \to \R^2$ with at most one cusp point
such that $f|_{S(f)}$ is a topological embedding.

In the case of $n$ odd, we will also show that
we can construct a $C^\infty$ stable map $f : M \to \R^2$
without cusps such that
the singular value set winds ``monotonically'' (with respect
to the polar angle) around the origin.
This leads to a new proof to the existence
of an open book structure on odd dimensional
smooth manifolds (see Theorem~\ref{thm:ob}),
which is originally due to Winkelnkemper
\cite{Wink}, Tamura \cite{Tamura}, Lawson
\cite{Lawson} and Quinn \cite{Quinn}.
In fact, we will show that any open book structure
on the boundary $\partial M$ of a compact odd dimensional manifold
$M$ extends through $M$.
The open book structure that we construct has a very
nice property with respect to a radial Morse function,
and this observation will be used to prove
Theorems~\ref{thm:ob2} and \ref{thm:ob3}:
if $M$ is a closed $k$--connected
$n$--dimensional manifold with $n \geq 7$ odd
and $1 \leq k \leq (n-5)/2$, then
there exists a $C^\infty$ stable map $f : M \to \R^2$
without cusp points such that 
$f$ has no fold points of absolute index $i$ with
$1 \leq i \leq k$.
This gives a reasonable answer to the second half
of Problem~\ref{prob1}.
It is surprizing to see such a generalization
of the above-mentioned result due to Smale, which
led to the solution to the higher dimensional
Poincar\'e Conjecture, to generic maps (or $C^\infty$
stable maps) into the plane,
although they are restricted to odd dimensions, since
now the singular points are not discrete and constitute
positive dimensional submanifolds.
The main ingredient for our construction of 
deformations of generic maps is the so-called
Cerf Theory \cite{Cerf, HW}.

The paper is organized as follows. In \S\ref{section2},
we review important properties of generic maps
of manifolds of dimension $n \geq 2$
into $2$--dimensional manifolds.
In particular, we give some lists of
changes of the singular value set for a generic
$1$--parameter family of such maps, which will
play an essential role in this paper.
In fact, in the case of manifolds $M$
of dimension $n = 4$, a similar set of lists
has been given in \cite{BS1} (see also \cite{BSPNAS}).

In \S\ref{section3}, we introduce the notion
of moves, which change the singular value set
of a generic smooth map. We rigorously formulate
the notion of an always-realizable move, which
was originally introduced in \cite{BS1},
as a formal change of the singular value set
that is always realized by a generic $1$--parameter
family of smooth maps.

In \S\ref{section4}, using such always-realizable moves,
we start to arrange the singular value set of generic
maps into $S^2$. We first arrange them so that
they wind monotonically around a certain point in the target,
and then we arrange them in the order of the
absolute index. These explicit procedures using
always-realizable moves lead to the proof of
Theorem~\ref{thm1}.

In \S\ref{section5}, using such a generic smooth
map into $S^2$, in the even dimensional case,
we construct a generic smooth map into $\R^2$
and prove Theorem~\ref{thm2}.
The idea of the proof of these theorems are
basically similar to that exploited in
dimension $4$ in \cite{BS1}. However, in dimension
$4$, one can use Lefschetz singularities
(see Remark~\ref{rem:Lef}). In our case
of arbitrary dimensions, we will give
new techniques that do not involve Lefschetz
singularities. This is one of the new contributions
of this paper.

In \S\ref{section6}, for a compact manifold $M$
of odd dimension with an open book structure
being given on the boundary,
we construct a smooth generic map into $\R^2$ 
which is compatible with the open book structure
on the boundary, and,
using techniques employed in \S\ref{section4},
we arrange the singular value set so that it
winds monotonically around the origin. Then, we prove
the existence of an open book structure on $M$ which
extends the given one on $\partial M$.

In \S\ref{section7}, we study the connectivity
of the page and the binding of the open book structure thus
constructed, and using the Cerf theory techniques
developed in \cite{HW}, we eliminate 
fold loci of certain absolute indices.

Finally, in \S\ref{section8}, we
consider the topological
position of the singular point set of a generic
map, instead of its image. Such study was
conducted in the author's previous paper \cite{Sa0};
however, in some cases, the result was not
correctly stated. In this section, we give
a correction to the result and its proof.

Throughout the paper, all manifolds and maps are
smooth of class $C^\infty$ unless otherwise specified.
The symbol ``$\cong$'' denotes a diffeomorphism
between manifolds or an appropriate isomorphism
between algebraic objects.
All homology groups are with integer coefficients unless
otherwise specified.
For a smooth map $f : M \to N$ between smooth manifolds,
$S(f) (\subset M)$ denotes the set of singular points of $f$.
For a real number $x$, $\lfloor x \rfloor$ denotes
the greatest integer not exceeding $x$.

\section{Generic maps into surfaces}\label{section2}

Let $M$ be a smooth closed manifold of dimension $n \geq 2$,
and $N$ a smooth surface.

\begin{defn}\label{def:fold-cusp}
A singular point $p$ of a smooth map $f : M \to N$
that can be described by the normal form
$$(t, x_1, x_2, \ldots, x_{n-1})
\mapsto (X, Y) =
(t, -x_1^2 - \cdots -x_i^2 + x_{i+1}^2 + \cdots + x_{n-1}^2)$$
with respect to appropriate coordinates around $p$ and $f(p)$
is called a
\emph{fold singularity} (or a \emph{fold}), where the integer
$i$ is called the \emph{index} with respect to the $(-Y)$--direction.
When we choose the $Y$--coordinate so that
$0 \leq i \leq \left\lfloor (n-1)/2 \right\rfloor$ holds,
$i$ is called the \emph{absolute index}.\footnote{We warn the
reader that this convention is different from that 
adopted in \cite[p.~273]{Lev}, \cite{Wr}, etc.}
Note that the absolute index of a fold is well-defined and
does not depend on the choice of the local coordinates.

When $i=0$ or $n-1$, the singular point is called a \emph{definite fold}.
Otherwise, it is called an \emph{indefinite fold}.
\end{defn}

\begin{defn}
A singular point of a smooth map $M \to N$
that has the normal form 
$$(t_1, t_2, x_1, x_2, \ldots, x_{n-2})$$
$$\mapsto (X, Y) =
(t_1, t_2^3 + t_1 t_2 -x_1^2 - \cdots -x_i^2 + x_{i+1}^2 + \cdots + x_{n-2}^2)$$
is called a
\emph{cusp singularity} (or a \emph{cusp}), 
where $i$ is called the \emph{index} with respect to the $(-Y)$--direction.
When we choose the $Y$--coordinate so that
$0 \leq i \leq \left\lfloor (n-2)/2 \right\rfloor$ holds,
$i$ is called the \emph{absolute index}, which is well-defined.
\end{defn}

A smooth map $f : M \to N$ that has only folds and cusps 
as its singularities
is called a \emph{generic map} (or an \emph{excellent map}). 
On the other hand, if $f$ has only fold singularities, then
it is called a \emph{fold map}.\footnote{In \cite{GG}, this is
called a \emph{submersion with folds}.}
It is well known that
a smooth map $M \to N$ can always be approximated by a generic map
arbitrarily well and hence is
homotopic to a generic map \cite{T, Wh}.

For a smooth map $f : M \to N$, 
$S(f)$ denotes the set of singular points of $f$.
If $f$ is generic and the source manifold $M$ is closed
(compact and without boundary), then $S(f)$
is a disjoint union of finitely many circles, consisting of
finitely many cusps together with finitely many arcs and circles of
folds, and
$f|_{S(f) \setminus \{\mbox{\footnotesize{cusps}}\}}$
is an immersion.

In the following, we assume that the source manifold
is closed unless otherwise stated.

\begin{defn}\label{def:stable}
A generic map $f : M \to N$ of an $n$--dimensional
manifold $M$ with $n \geq 2$ into a surface $N$
is a \emph{$C^\infty$ stable map} (or a \emph{stable
map}, for short) if the following conditions
are satisfied.
\begin{enumerate}
\item The restriction $f|_{S(f) \setminus \{\mbox{\footnotesize{cusps}}\}}$
is an immersion with normal crossings, and
\item for every cusp $p \in M$, we have $f^{-1}(f(p)) \cap S(f) = \{p\}$.
\end{enumerate}
\end{defn}

Note that $C^\infty$ stable maps are usually defined
in another way (see \cite{GG}): however, it is known that
the above definition is equivalent to such a usual
definition (see, for example, \cite{Lev1}).
Note also that the set of stable maps is open and dense in the
space $C^\infty(M, N)$ of all smooth maps endowed with the Whitney
$C^\infty$--topology.

\begin{figure}[h]
\centering
\psfrag{S}{$f(S(f))$}
\psfrag{L}{$\ell$}
\psfrag{N}{$N$}
\psfrag{f}{$f(p)$}
\includegraphics[width=0.9\linewidth,height=0.2\textheight,
keepaspectratio]{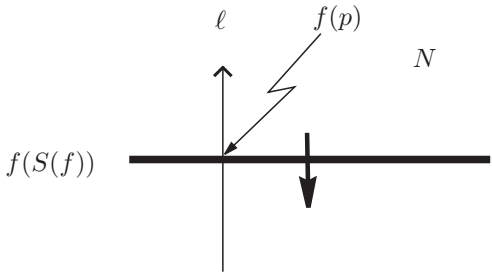}
\caption{Normal orientation for $f|_{S(f)}$, depicted
as thick downward arrow}
\label{fig602}
\end{figure}

Let $f : M \to N$ be a generic map.
For a fold $p \in M$ of index $i$ with respect
to the $(-Y)$--direction in the sense of Definition~\ref{def:fold-cusp},
the immersion $f|_{S(f)}$ near $p$
is \emph{normally oriented} so that $f|_{f^{-1}(\ell)} : 
f^{-1}(\ell) \to \ell$
has $p$ as a non-degenerate critical point of index $i$, where
$\ell$ is a small open arc embedded in $N$ which intersects
$f(S(f))$ transversely at exactly one point $f(p)$, and
the normal orientation corresponds to the negative direction of $\ell
\cong \R$, i.e.\ the $(-Y)$--direction (see Fig.~\ref{fig602}).

Note that if we reverse the normal orientation, then the
index $i$ changes to $(n-1)-i$. Therefore, if we require that
$i$ should be
the absolute index, then the normal orientation
is uniquely determined, except for the case where $n$ is odd and
$i = (n-1)/2$.
When the absolute index satisfies $i = (n-1)/2$,
both of the normal orientations give the same index.
In this case, we usually do not associate any normal orientation
(or, in other words, both normal orientations are associated, or
the normal orientation is redundant).

\begin{figure}[h]
\centering
\psfrag{a}{$i$}
\psfrag{b}{$i+1$}
\psfrag{c}{$i$}
\includegraphics[width=0.9\linewidth,height=0.15\textheight,
keepaspectratio]{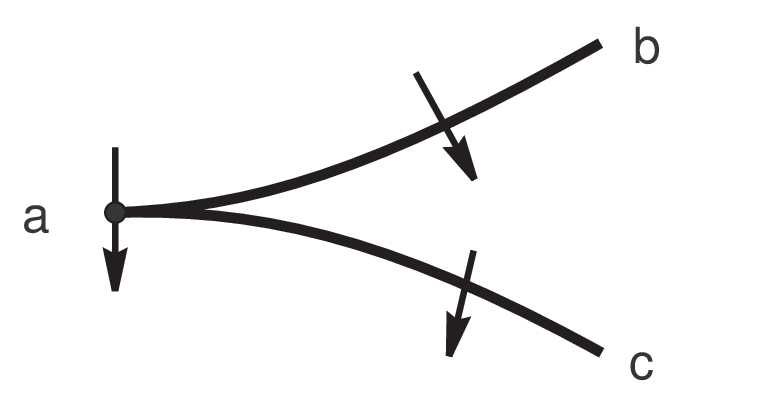}
\caption{Indices and normal orientations around the image of a cusp}
\label{fig614}
\end{figure}

For the image of a cusp point, we can also
consider a ``normal orientation'': this
is considered to be normal to the two branches
of curves emanating from the point. The indices of the folds
around the image
of a cusp are always as depicted
in Fig.~\ref{fig614} for some integer $i$ with respect to
appropriate normal orientations as in the figure.

\begin{figure}[h]
\centering
\psfrag{a}{$i$}
\psfrag{b}{$i+1$}
\psfrag{c}{$i$}
\psfrag{a2}{$(n-3)/2$}
\psfrag{b2}{$(n-1)/2$}
\psfrag{c2}{$(n-3)/2$}
\psfrag{a3}{$(n-2)/2$}
\psfrag{b3}{$(n-2)/2$}
\psfrag{c3}{$(n-2)/2$}
\psfrag{1}{(1) $i + 1 < (n-1)/2$}
\psfrag{2}{(2) $i + 1 = (n-1)/2$}
\psfrag{3}{(3) $i = (n-2)/2$}
\includegraphics[width=0.9\linewidth,height=0.5\textheight,
keepaspectratio]{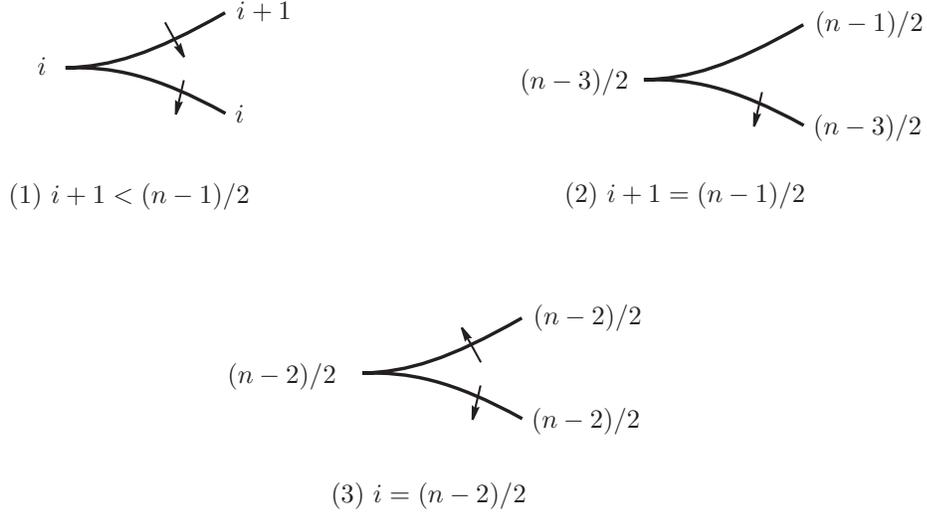}
\caption{\emph{Absolute indices} and normal orientations around cusp images}
\label{fig601}
\end{figure}

We see easily that the absolute indices around cusps and the
corresponding normal orientations are as indicated in
Fig.~\ref{fig601}. Note that when $n$ is odd, Case~(2) can occur, 
while when $n$ is even, 
Case~(3) can occur. Here, we usually do not
indicate the normal orientation of the image
of a cusp. This is because for Case (1), the normal
orientation is considered to be the unique one
consistent with the normal orientations for the
two branches of curves (as in Fig.~\ref{fig614}),
for Case (2), it should be considered to be
the unique one consistent with the normal
orientation of one of the two branches, and 
for Case (3), both of the normal orientations can
be chosen.

\section{Diagrams and their moves}\label{section3}

In this section, we introduce the notion of 
diagrams and moves for them.
In the following, we assume that the source manifold
$M$ is closed of dimension $n \geq 2$ and that $N$ is a smooth surface.

\begin{defn}\label{defn:diagram}
Let $S$ be a finite disjoint union of circles
and $\gamma : S \to N$ a smooth map with
the following properties.
\begin{enumerate}
\item The map $\gamma$ is an immersion with normal crossings
outside of a finite set $C (\subset S)$.
\item Around each point $c_\ell$ of $C$, $\gamma$ is of the form $t \mapsto (t^2, t^3)$
with respect to appropriate local coordinates around $c_\ell$ and $\gamma(c_\ell)$.
\item For each point $c_\ell \in C$, we have $\gamma^{-1}(\gamma(c_\ell)) = 
\{c_\ell\}$.
\end{enumerate}
We also assume that the immersion obtained by
restricting $\gamma$ to each component $S_k$ of $S \setminus C$
is normally oriented. 
Moreover, to each $S_k \subset S \setminus C$ is associated
an integer $i_k$ with $0 \leq i_k \leq \lfloor (n-1)/2 \rfloor$, and
to each $c_\ell \in C$ is associated
an integer $j_\ell$ with $0 \leq j_\ell \leq \lfloor (n-2)/2 \rfloor$
in such a way that around each point of $C$, they are
as depicted in Fig.~\ref{fig601}. When $n$ is odd and $S_k$ has
index $i_k = (n-1)/2$, it is allowed that no normal orientation
is given for $\gamma|_{S_k}$.
Then, we denote by $\mathcal{D}$
the image $\gamma(S)$ together with
the associated normal orientations and the integers as above,
and we call the pair $(N, \mathcal{D})$
(or simply, $\mathcal{D}$) a
\emph{normally oriented indexed diagram} (or a \emph{diagram}, in short).
For each point $c_\ell \in C$, we call its image $\gamma(c_\ell)$ a 
\emph{cuspidal point} of $\mathcal{D}$.

For example, if $f : M \to N$ is a stable map, then
$(N, f(S(f)))$ together with the normal orientations
and the absolute indices as described in the previous section
can be considered to be a diagram.

In order to make the descriptions simpler, we adopt
the convention that 
\begin{enumerate}
\item for $S_k$, the associated integer $i_k$ can be chosen
from the set $\{0, 1, \ldots, n-1\}$,
\item for $c_\ell$, the associated integer $j_\ell$
can be chosen from the set $\{0, 1, \ldots, n-2\}$, 
\item if we reverse the normal orientation of $\gamma|_{S_k}$,
then the index $i_k$ turns into $(n-1)-i_k$,
and the indices $j_\ell$ for the adjacent cuspidal points turn
into $(n-2)-j_\ell$:
i.e., such
two diagrams are considered to be the same.
\end{enumerate}
Note that for a cusp point, we often omit the index and we do
not associate a normal orientation.
This is because the normal orientation and 
the index of a cuspidal point can be read-off from the
normal orientations and the indices of the two
adjacent branches of curves. 

Furthermore, two normally oriented indexed diagrams are
\emph{isotopic} if there exists an ambient isotopy of $N$
which takes one to the other preserving the normal orientations
and the indices. We often use the terminology ``diagram''
in order to refer to its isotopy class.
\end{defn}

\begin{defn}\label{defn:move}
A \emph{move} for a normally oriented indexed
diagram $(N, \mathcal{D})$ is an operation
which changes $\mathcal{D}$ to a new diagram
$\mathcal{D}'$ such that
\begin{enumerate}
\item there exists a small disk $\Delta$ embedded
in $N$ such that $\partial \Delta$ intersects $\mathcal{D}$
and $\mathcal{D}'$ transversely outside of their cuspidal points
and double points, and
\item $\mathcal{D}$ coincides with $\mathcal{D}'$
outside of $\Delta$, together with the normal orientations
and the indices.
\end{enumerate}
Since such an operation does not change anything
outside of the disk $\Delta$,
the general operation which replaces
$\Delta \cap \mathcal{D}$ with $\Delta \cap \mathcal{D}'$,
is also called a \emph{move}, and it is often represented
by the symbol ``$\Delta \cap \mathcal{D} \longrightarrow \Delta \cap \mathcal{D}'$.''
\end{defn}

\begin{defn}
Suppose that a move $\Delta \cap \mathcal{D} \longrightarrow 
\Delta \cap \mathcal{D}'$
is given. Let $f : M \to N$ be an arbitrary stable
map such that 
\begin{itemize}
\item[$(\sharp)$] there exists a disk $\Delta_f$ in $N$
such that for the diagram $(N, f(S(f)))$, 
$\partial \Delta_f$ intersects $f(S(f))$ outside of
the images of the cusps and double points transversely and that
under a certain identification of $\Delta_f$ with $\Delta$,
$\Delta_f \cap f(S(f))$ is isotopic to $\Delta \cap \mathcal{D}$,
relative to the boundary $\partial \Delta_f \cong \partial \Delta$.
\end{itemize}
Suppose that there exists a smooth
$1$--parameter family $f_t : M \to N$, $t \in [0, 1]$, 
of smooth maps, with $f_0 = f$, which modifies the singular point
set image $f(S(f))$ only in $\Int{\Delta_f}$,
such that $\Delta_f \cap f_1(S(f_1))$ is isotopic to
$\Delta \cap \mathcal{D}'$ relative to the boundary $\partial \Delta_f
\cong \partial \Delta$,
that $f_t^{-1}(\Delta_f) = f^{-1}(\Delta_f)$ for all $t$, that
$f_t=f$ outside of the inverse images of $\Delta_f$ for all $t$,
and that except for finitely many $t$'s, $f_t$ is a stable map.
Then, we say that the move
$\Delta \cap \mathcal{D} \longrightarrow \Delta \cap \mathcal{D}'$
is \emph{realized} for $f$ (or is \emph{realizable}).
If a move is realized for every $f$ with property $(\sharp)$ above,
then we say that the move is \emph{always-realizable}.
\end{defn}

In the following, we give several sets of moves that
appear in a generic $1$--parameter family of
smooth maps of closed $n$--dimensional
manifolds, $n \geq 2$, into surfaces.
Here, a smooth $1$--parameter family of smooth
maps $\{f_t\}$ is said to be \emph{generic} if
every $f_t$ is a stable map except
for a finite number of $t$'s.
We usually impose a more technical condition
that the family corresponds to a
curve in the mapping space that is transverse to 
codimension $1$ unstable loci.
For more technical details, see \cite{Ch1, Ch2, GK2, La, Soto}, 
for example.

Some of the following moves are always-realizable, while others are not.
The reader should refer to the arguments given
in \cite{BS1}, where the case of $n=4$ is thoroughly explained.
The arguments for general dimension $n$ are basically similar,
although there are some subtle issues for higher dimensions.
For details, the reader is referred to \cite{Cerf, CL, HW}.

\begin{figure}[h]
\centering
\psfrag{E}{$\emptyset$}
\psfrag{B}{Birth}
\psfrag{D}{\bang Death}
\psfrag{FM}{\bang Fold Merge}
\psfrag{CM}{Cusp Merge}
\psfrag{CM1}{(If fibers over $(\ast)$ are connected.)}
\psfrag{F}{Flip}
\psfrag{UF}{Unflip}
\psfrag{UF1}{(If $i=0$ or $2 \leq i \leq n-2$.)}
\psfrag{i}{$i$}
\psfrag{i1}{$i+1$}
\psfrag{st}{$(\ast)$}
\includegraphics[width=\linewidth,height=0.6\textheight,
keepaspectratio]{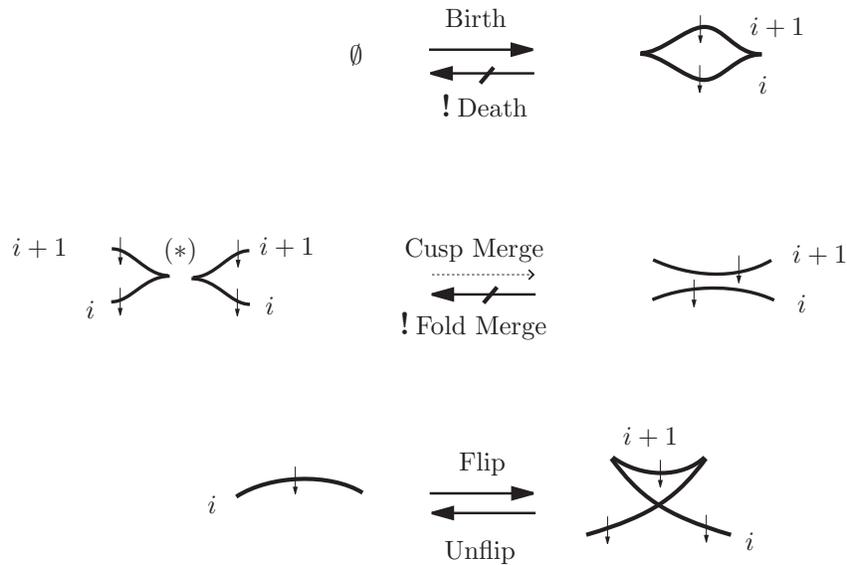}
\caption{Moves of type I}
\label{fig604}
\end{figure}

Fig.~\ref{fig604} gives the first set of moves that
we will use later. Note that
outside of the small disks depicted in the figure,
these moves do not change the diagrams.
These are called \emph{moves of type I}: note that each move
corresponds to a generic $1$--parameter family of
smooth maps that
changes the map only in a neighborhood of a point
in the source manifold.
The reader is referred to \cite[Part~I, Chapter~V and 
Chapter~VII, Proposition~1]{HW}
for $n \geq 6$.

Note that the simple arrows are always-realizable, while
the arrows with a slanted line are not.
``Cusp Merge'', which is depicted by a dotted arrow, 
is ``always-realizable'' as long as the fibers over the
region indicated by $(\ast)$ of a given stable map
are connected. Furthermore, ``Unflip'' is `always-realizable''
as long as the index satisfies the assumption
$i=0$ or $2 \leq i \leq n-2$ (see 
\cite[Part~I, Chapter~V, Proposition~1.4]{HW}).

\begin{figure}[h]
\centering
\psfrag{i}{$i$}
\psfrag{j}{$j$}
\psfrag{ij}{$i < j$}
\psfrag{ji}{$i > j$}
\psfrag{iijj}{$i \leq j$}
\includegraphics[width=\linewidth,height=0.4\textheight,
keepaspectratio]{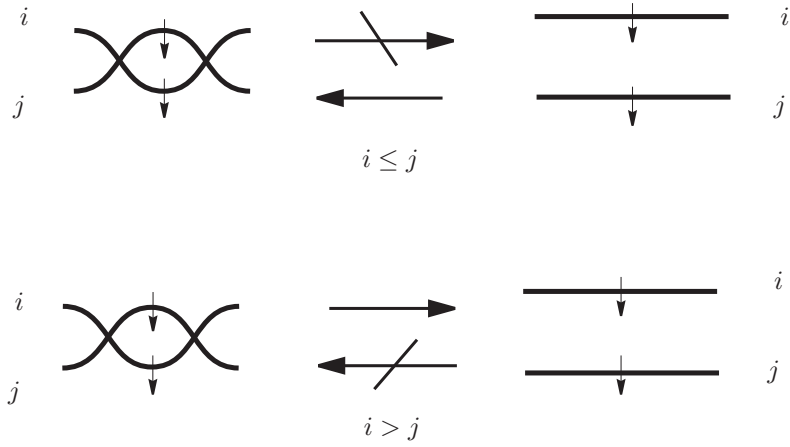}
\caption{Moves of type II}
\label{fig605}
\end{figure}

The moves as depicted
in Fig.~\ref{fig605} are called \emph{moves of type II}.
Note that these moves correspond to generic $1$--parameter
families of smooth maps that change the maps
in neighborhoods of
two points in the source manifold.
For these moves and the following, the reader is referred to
\cite[Part~I, Chapter~V and Chapter~I, \S\S7 and 8]{HW} and 
\cite[Chapitre~IV, \S2.2, Proposition~2]{Cerf}.

The important feature of the moves of type II (and those
of type III that follows)
is that fold image of smaller index can be pushed
to the direction of the normal orientation.

\begin{figure}[t]
\centering
\psfrag{i}{$i$}
\psfrag{j}{$j$}
\psfrag{k}{$k$}
\psfrag{jik}{$i+k \geq n$, $\max\{i,k\} > j$}
\psfrag{jik1}{or $i=j=k \geq 2$}
\includegraphics[width=0.7\linewidth,height=0.3\textheight,
keepaspectratio]{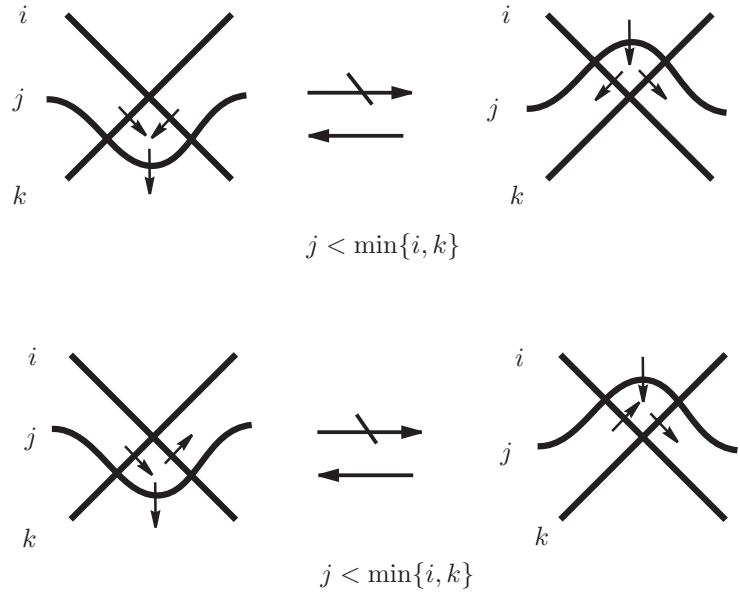}
\caption{Moves of type III}
\label{fig606}
\end{figure}

The moves as depicted in Fig.~\ref{fig606}
are called \emph{moves of type III}.
Note that these moves correspond to
generic $1$--parameter families of smooth maps
that change the maps in neighborhoods of
three points in the source manifold.

\begin{figure}[t]
\centering
\psfrag{i}{$i+1$}
\psfrag{j}{$i$}
\psfrag{k}{$j$}
\psfrag{k0}{$j > 0$}
\psfrag{k1}{$j < i$ or}
\psfrag{k2}{$j < n-2-i$}
\includegraphics[width=\linewidth,height=0.5\textheight,
keepaspectratio]{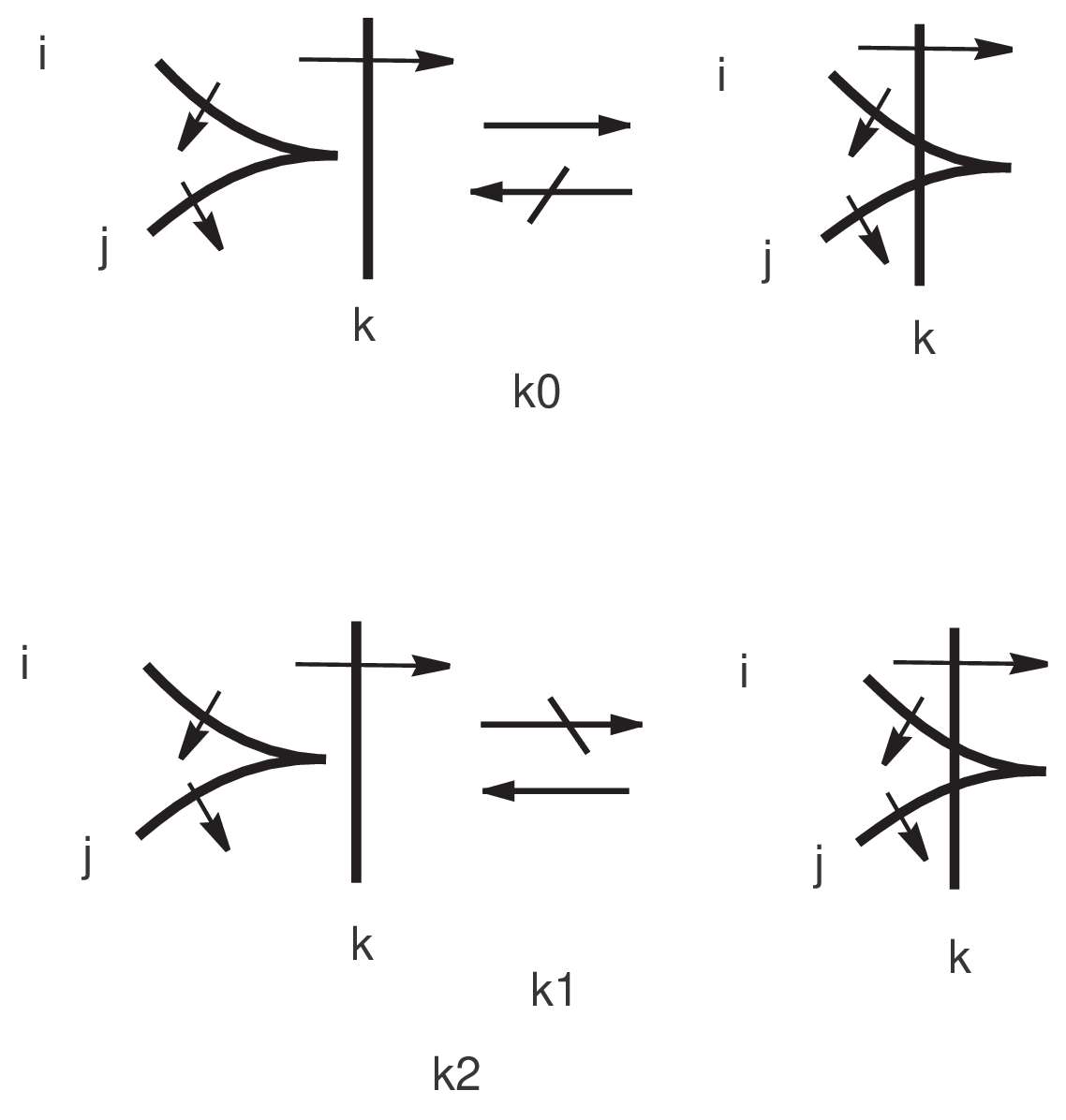}
\caption{Cusp--fold crossing moves}
\label{fig607}
\end{figure}

The moves as depicted in Fig.~\ref{fig607} are called
\emph{cusp--fold crossing moves} (refer to \cite[Part~I,
Chapter~V, (0.3)]{HW}
and \cite[Chapitre~IV, Proposition~4]{Cerf}).

\begin{figure}[t]
\centering
\psfrag{II}{exchange}
\psfrag{II1}{(If fibers over $(\ast)$ are connected.)}
\psfrag{a}{$(\ast)$}
\psfrag{cc}{criss-cross}
\psfrag{cc1}{(If fibers over $(\ast)$ are connected.)}
\psfrag{i}{$i$}
\includegraphics[width=\linewidth,height=0.5\textheight,
keepaspectratio]{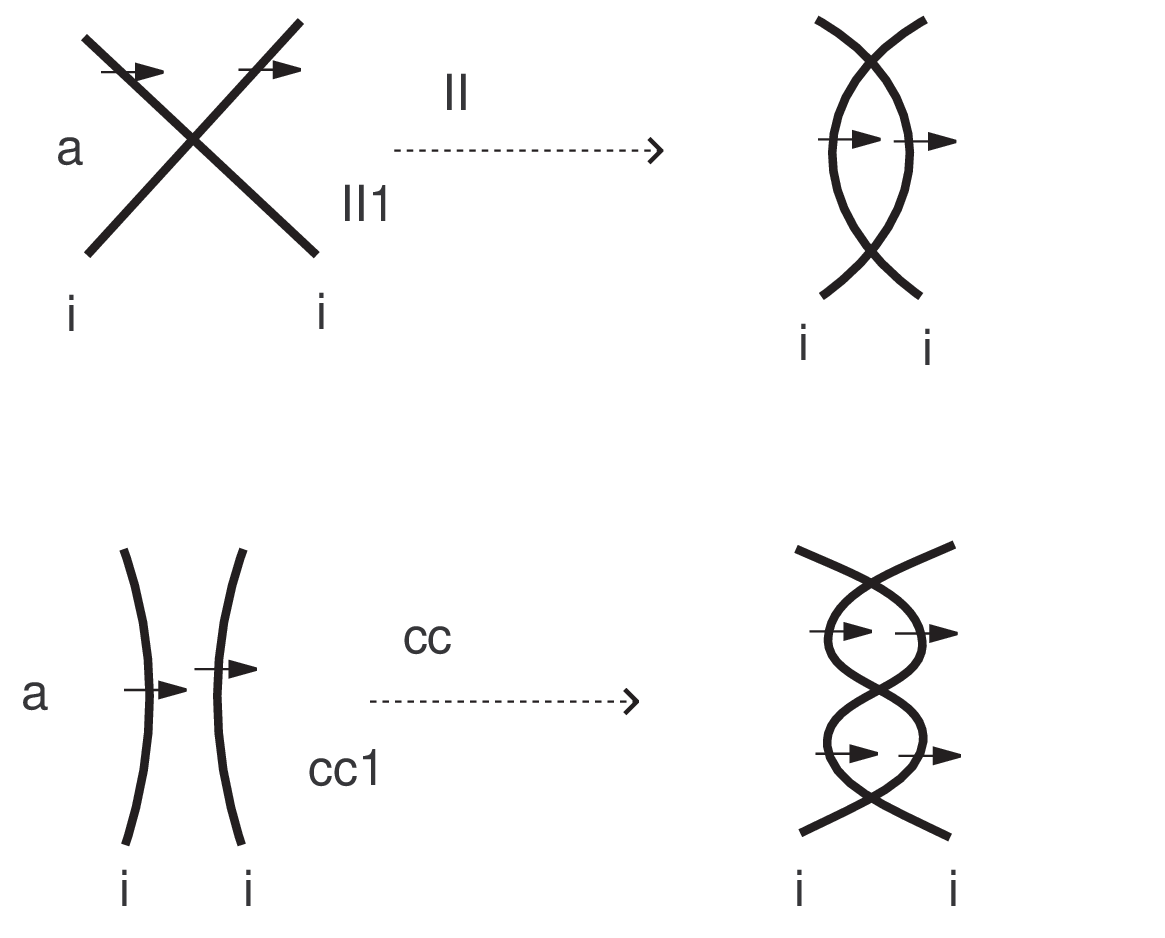}
\caption{Exchange move and criss-cross move, which are valid for $i \geq 1$.}
\label{fig609}
\end{figure}

For our purpose,
we also need some combination moves which
are realized by a finite sequence of realizable moves.
They are called an \emph{exchange move} and a
\emph{criss-cross move} (see Fig.~\ref{fig609}).
Note that these moves are
``always-realizable'' as long as the fibers are connected
over the regions indicated by $(\ast)$ of a given stable map.
Therefore, they are depicted by dotted arrows. 

\begin{figure}[t]
\centering
\psfrag{fi}{flips}
\psfrag{c}{cusp}
\psfrag{m}{merge}
\psfrag{III}{$\mathrm{III}$}
\psfrag{f}{unflip}
\psfrag{i}{$i$}
\includegraphics[width=\linewidth,height=0.25\textheight,
keepaspectratio]{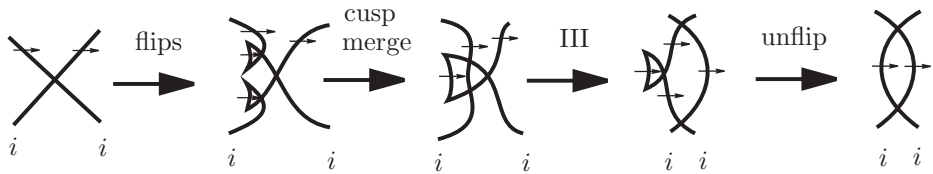}
\caption{Moves for the exchange move, $i \geq 1$.}
\label{fig610}
\end{figure}

An exchange move is realized as in Fig.~\ref{fig610}.
A criss-cross move is realized by combining
a move of type II and an exchange move.

A more detailed discussion about
the exchange move can be found in
\cite[Proposition~3.2]{BS1}
for $n=4$ and $i=1$. 
For higher dimensions, the case of $i=1$
can be proved by the same argument as
in \cite{BS1}. For $i \geq 2$, we can apply
flips, cusp merges, a move of type III, and
an unflip safely, since in our situation they
are all ``always-realizable.''

\section{Stable maps into $S^2$}\label{section4}

In this section, we state and prove our first
main result concerning stable maps into $S^2$.
Let us start with the following definition.

\begin{defn}
A generic map $f : M \to N$ of a closed
manifold into a surface is said to be
\emph{image simple}
if $f|_{S(f)}$ is a topological embedding.
\end{defn}

Note that an image simple generic map
is always stable as long as the source manifold
is closed.

\begin{rk}
Recall that a generic map $f : M \to N$ is
\emph{simple} if it has no cusps and 
each connected component
of every fiber contains at most one singular point
(see \cite{Lev1}).
In the literature, some authors use the terminology
``\emph{strongly simple}'' instead of ``image simple''.
Note that a generic map which is image simple and
has no cusps is simple.
\end{rk}

In this section, we prove the following.

\begin{thm}\label{thm1}
Let $M$ be a closed connected $n$--dimensional manifold with $n \geq 3$.
Then, every map $M \to S^2$ is homotopic to a stable
map $f : M \to S^2$ with the following properties.
\begin{enumerate}
\item When $n$ is odd, $f$ has only folds of 
absolute index $(n-1)/2$ and has no cusps.
\item When $n$ is even, $f$ has only folds of absolute
index $(n-2)/2$ and at most one cusp, whose absolute index is equal to
$(n-2)/2$; furthermore, $f$ is image simple.
\end{enumerate}
\end{thm}

\begin{rk}\label{remarkn2}
Theorem~\ref{thm1} holds also for $n = 2$.
See \cite[Theorem~1.4]{TY}, for example.
\end{rk}

\begin{rk}
As the following proof shows, starting
from an arbitrary stable map $M \to S^2$,
we have an explicit algorithm, as a sequence
of always-realizable moves, that modifies
the given map to the one as described in the theorem.
\end{rk}

For the moment, let us consider
a stable map $f : M \to N$, where $M$
is a closed manifold of dimension $n \geq 3$
and $N = S^2$ or $\R^2$.
In the following, for a closed manifold $M$, $\chi(M)$ denotes the 
Euler characteristic of $M$.

According to Thom \cite{T}, the number of cusp points
of $f : M \to N$ has the same parity as $\chi(M)$.
Levine \cite{Lev} showed that if $M$ is connected, then
every map $M \to N$ is homotopic
to a generic map $f$ such that
\renewcommand{\theenumi}{\roman{enumi}}
\begin{enumerate}
\item $f$ has no cusp points when $\chi(M)$ is even, and
\label{c1}
\item $f$ has exactly one cusp point, whose index is equal to $(n-2)/2$, 
when $\chi(M)$ is odd (and in this case $n$ is necessarily even).
\label{c2}
\end{enumerate}
\renewcommand{\theenumi}{\roman{arabic}}

In fact, if we start with a stable map, then
Levine showed that by applying always-realizable
cusp-fold crossing moves and cusp merge
moves, one
can modify any given stable map homotopically to get one satisfying
the above conditions (\ref{c1}) and (\ref{c2}).

Therefore, for the proof of Theorem~\ref{thm1},
we may assume
that $f$ has at most one cusp of index $(n-2)/2$
when $n$ is even, and no cusp when $n$ is odd.

\begin{figure}[htbp]
\centering
\psfrag{f}{flip}
\psfrag{u}{unflip}
\psfrag{r}{II}
\psfrag{i}{$i_0$}
\psfrag{i1}{$i_0+1$}
\includegraphics[width=\linewidth,height=0.5\textheight,
keepaspectratio]{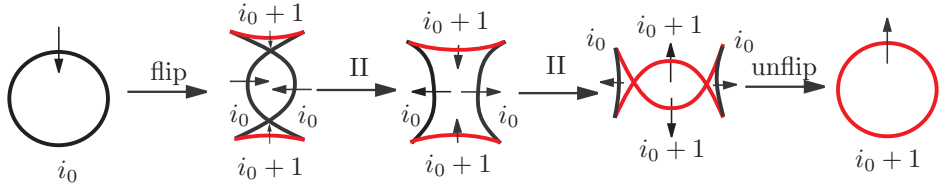}
\caption{Turning a fold circle of index $i_0$ to one of index $i_0+1$, which
is possible when $i_0 \leq (n-3)/2$.}
\label{fig611}
\end{figure}

Let us now consider the
elimination of definite folds.
The author \cite{Sa, Sa2} showed that
for a closed manifold $M$ of dimension $n \geq 3$, 
every map $M \to S^2$ is homotopic to
a stable map without definite fold singularities.
In fact, the proof shows that any given stable
map with at most one cusp point
can be modified by applying 
some flips, cusp merges, 
moves of types II and III,
and the moves as depicted in Fig.~\ref{fig611} with
$i_0=0$ so that we get a stable map with at most
one cusp point and without definite folds.

Note that, clearly, this is, in general, impossible for maps into $\R^2$.

\begin{rk}
Let us consider the neighborhood of the index $0$
fold locus corresponding to Fig.~\ref{fig611} with $i_0=0$.
If it is orientable, then we see that the newly created fiber
over a point in the innermost region is diffeomorphic to
$S^1 \times S^{n-3}$. If it is non-orientable, then
we see that the fiber is diffeomorphic to the non-orientable
$S^{n-3}$--bundle over $S^1$.
\end{rk}

Then, for $n=3$, we see that the proof of Theorem~\ref{thm1}
is complete. In the following, we assume $n \geq 4$.

Let us assume that $f : M \to S^2$ is a stable
map with at most one cusp such that the minimum $i_0$
of the absolute indices of folds satisfies
$0 < i_0 < \lfloor (n-1)/2 \rfloor$.
We successively eliminate fold singularities (by homotopy)
in the order of absolute index, without changing the number
of cusps, as follows.

First, let us define the following notion.

\begin{defn}
Let $f: M \to N$ be a stable map into a surface $N$.
We say that the fold image of $f$ of absolute index $i$
is \emph{outward-directed} (resp.\ \emph{inward-directed}), 
if the image of all the folds of absolute index $i$ of $f$
is contained in the interior of
a $2$--disk $D$ in $N$ in such a way that the complement of 
a regular value $z_0 \in D$ can be nonsingularly
foliated by arcs oriented from $z_0$ to $\partial D$, 
which intersect the image of each
fold arc of absolute index $i$ transversely in its normal direction 
(resp.\ the opposite direction).
In this case, the connected component of the
complement of the fold image of absolute index $i$
in $D$ containing $z_0$ is called the
\emph{innermost region}.
\end{defn}

Note that when $f$ has exactly one cusp point of
index $(n-2)/2$, if the fold image of absolute
index $(n-2)/2$ of $f$ is outward-directed, then
the cusp image also points ``outward'' 
(see Fig.~\ref{fig601} (3)).

We first show the following, which makes the fold image
of absolute index $i_0$ outward-directed. We warn the reader that
the case with $i_0 = \lfloor (n-1)/2 \rfloor$ is excluded.

\begin{lem}\label{lem:directed}
Let $f : M \to S^2$ be a stable map
with at most one cusp point such that
the minimum $i_0$ of the absolute indices of folds
satisfies $0 < i_0 < \lfloor (n-1)/2 \rfloor$.
Then, by using always-realizable moves
we can make the fold image of absolute index $i_0$
outward-directed in such a way that the other
folds have absolute index strictly greater
than $i_0$, that their images are situated in the innermost region
and that the resulting stable map has at most
one cusp.
\end{lem}

\begin{proof}
Let $D \subset S^2$ be a $2$--disk containing
$f(S(f))$. We can also take an annulus $A
\cong [0, 1] \times S^1$ in $D$ containing $f(S(f))$
in its interior.
We use almost the same argument
as in the proof of \cite[Theorem~4.1]{BS1}
in order to deform the given stable map
so that the fold image of absolute
index $i_0$ is contained in $A$, and is outward-directed, and that
the other fold images are 
situated in the inner region of $D \setminus A$.
The points different from the situation
in \cite{BS1} that we need to take care of are
as follows.

(1) In \cite{BS1}, the dimension of the source
manifold $n$ was equal to $4$ and we had
only folds of absolute index $i_0=1$. In our present
case, the dimension is greater than or equal to $4$
and we may have folds of absolute index $> i_0$.
When we apply always-realizable moves as in \cite{BS1},
as the folds other than those of absolute index $i_0$
have absolute index strictly greater than $i_0$,
we can safely apply the necessary always-realizable moves.

(2) In \cite{BS1}, there were Lefschetz singular
points and one could use such singular points
in order to eliminate cusps, while in our situation we 
have no Lefschetz singular points, but may have
a cusp point. If we have a cusp, then $n$ is even and
the absolute index of the cusp
must be equal to $(n-2)/2$, which is maximal.
Therefore, we can use cusp-fold crossing
moves that are always-realizable instead of 
always-realizable moves of type II.

(3) In \cite{BS1}, after applying a flip that creates
two cusps, we eliminated the cusps by using
Lefschetz singular points. In our situation,
we can realize the same procedure as in Fig.~\ref{fig613}
by creating a fold circle of absolute index $i_0+1$
whose image is embedded and is outward-directed.
Note that the cusp merge can be performed, since
there exists a small $2$--disk containing the two
cusp images such that the two cusps lie on the
same connected component
of the inverse image of the $2$--disk.

\begin{figure}[h!]
\centering
\psfrag{c}{cusp merge}
\psfrag{f}{flip}
\psfrag{i}{$i_0$}
\psfrag{i1}{$i_0+1$}
\psfrag{2}{II}
\includegraphics[width=\linewidth,height=0.4\textheight,
keepaspectratio]{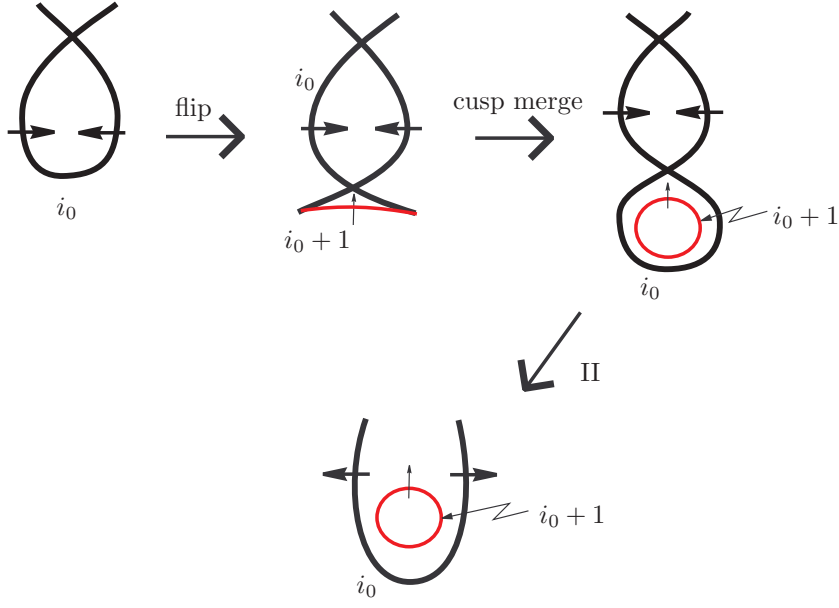}
\caption{Eliminating a self fold crossing}
\label{fig613}
\end{figure}

Details of the rest of the proof are left to the reader.
\end{proof}

\begin{rk}\label{rk:directed}
(1) In fact, we could perform the above procedure
before eliminating the definite folds. In that case, we
get a stable map as in the lemma with $i_0 = 0$.

(2) The above lemma is for maps into $S^2$: however,
as the proof shows, the necessary modifications are performed
only in the inverse images of the $2$--disk $D$.
The map outside that part is unchanged. Therefore, for example,
the regular fibers outside of $D$ stay the same.

(3) When $i_0 = \lfloor (n-1)/2 \rfloor$, the arguments
in the proof does not work, since $i_0+1$ is
no longer an absolute index. Furthermore, when $n$ is odd
and $i_0 = \lfloor (n-1)/2 \rfloor$,
we cannot apply type II moves as in the above
argument, since reversing the normal orientation
of the index $i_0$ fold image
does not change the index.
\end{rk}

Then, we need the following in order to eliminate
folds of absolute index $i_0$.

\begin{lem}\label{lem:elim}
Let $f : M \to S^2$ be a stable map
with at most one cusp point such that
the minimum $i_0$ of the absolute indices of folds
satisfies $0 < i_0 < \lfloor (n-1)/2 \rfloor$.
We assume that $M$ is connected.
Then, by using always-realizable moves
we can eliminate the folds of absolute index $i_0$
in such a way that the other
folds have absolute index strictly greater
than $i_0$ and that the resulting stable map has at most
one cusp.
\end{lem}

\begin{proof}
First, by Lemma~\ref{lem:directed},
we may assume that 
the fold image of absolute index $i_0$
is outward-directed in such a way that the other
folds have absolute index strictly greater
than $i_0$, and that their images are situated in the innermost region,
which is homeomorphic to an open disk and is denoted by $R$.
Since $f$ has directed fold image of index $i_0$, 
we may assume that it winds around the north pole $NP \notin R$
of $S^2$ and that it is normally oriented toward $NP$.

Note that if there is a cusp, then the unique cusp
lies on a component adjacent to fold loci whose
absolute indices are strictly greater than $i_0$.
In the following argument, we only modify the fold
image of index $i_0$, so that we do not move
the fold image adjacent to the cusp image.

Now, observe that the fibers over the points in $R$
are all connected.
This can be seen as follows. As the absolute indices
of folds whose images lie inside $R$ are greater
than or equal to $i_0+1$, we see that when we cross
a fold image, the topological change of fibers
can be described by a surgery corresponding to
the attachment of a handle
of index $i$ with $2 \leq i \leq n-3$. Therefore, the number of
connected components of fibers does not change. 
If we consider the topological changes of fibers
toward $NP$, then they correspond to 
the attachments of handles of index $n-1-i_0$, since the fold image of
index $i_0$ is normally oriented toward $NP$.
This means that the numbers of connected components
of fibers may increase (if $i_0 = 1$); however,
they never strictly decrease. 
Therefore, if a fiber over a point in $R$ is
disconnected, then $f^{-1}(R)$
is disconnected, since $R$ is simply connected;
hence, the manifold $M$
must necessarily be disconnected.
This is a contradiction.
Therefore, the fibers over the points in $R$
are all connected.

Then, by an argument similar to that in \cite{BS1},
we can arrange the fold image of index $i_0$ so that the
image of each fold component
of index $i_0$ winds exactly once
around $NP$, by using moves of type II together with criss-cross
moves.

We further follow the procedure as described 
in \cite[Figs.~30 and 31]{BS1}.
If the top strand, denoted by $b$, 
does not intersect the other strands, then we
can isotope it into a small open disk neighborhood of the north pole $NP$.
Then, we can use the move as described in Fig.~\ref{fig611}
of the present paper so as to turn it to a fold
image of absolute index $i_0+1$ which is an embedded
small circle and is directed outward with respect to $NP$.
If the top strand does intersect the others,
then we perform
the moves as described in \cite[Fig.~33]{BS1},
by using the move as in Fig.~\ref{fig613}
of the present paper instead of the move
as described in \cite[Fig.~34]{BS1}.
Note that by this move, we again create an embedded fold
image circle of index $i_0+1$, which can
be isotoped to a neighborhood of the south pole $SP \in R$.
In order to reach the second stage of
\cite[Fig.~33]{BS1}, we need to swing
the relevant fold image of index $i_0$ over
the region around $SP$: this is
possible, since all the fold images
in this region have absolute index 
strictly greater than $i_0$.

Then, we reach the final stage of \cite[Fig.~35]{BS1}.
Finally, we apply the move as in Fig.~\ref{fig611}
of the present paper to turn the index $i_0$ circle
to one of index $i_0+1$.

Repeating these procedures for the rest
of the fold components of absolute index $i_0$,
we can eliminate all the folds of absolute
index $i_0$ so that all the folds
of the resulting stable
map have absolute index strictly greater than $i_0$
and that it has at most one cusp.
This completes the proof.
\end{proof}

When the dimension of the source manifold is even, we can
further arrange the images of the folds and cusps as follows.

\begin{lem}\label{lem:emb}
Suppose that $n \geq 4$ is even.
Let $f : M \to S^2$ be a stable map
with at most one cusp point such that
the absolute indices of all folds
are equal to $i_0 = \lfloor (n-1)/2 \rfloor = (n-2)/2$.
We assume that $M$ is connected.
Then, by using always-realizable moves
we can make the fold image topologically embedded
in such a way that all the folds have absolute
indices equal to $i_0$ and that 
the resulting stable map has at most
one cusp.
\end{lem}

\begin{proof}
We follow the same argument as in the proof of
Lemma~\ref{lem:elim}. However, we need to
be careful, since if there is a cusp, then
we need to move the fold loci adjacent to
the unique cusp as well.

First, let us consider the case $n \geq 6$.
If we proceed as in the proof of Lemma~\ref{lem:elim}, then
there may appear some small embedded fold image circles
of index $i_0+1$ normally oriented outward
(see Figure~\ref{fig613}).
Note that if we reverse the normal orientation,
then the index turns into $i_0$, the absolute index,
since $n$ is even and $i_0 = (n-2)/2$.
We used push moves under the presence of Lefschetz
critical points (see \cite[Fig.~2]{BS1}).
In our present case, we can instead
use the moves as depicted in Fig.~\ref{fig617}.

\begin{figure}[h!]
\centering
\psfrag{i}{$i_0$}
\includegraphics[width=\linewidth,height=0.4\textheight,
keepaspectratio]{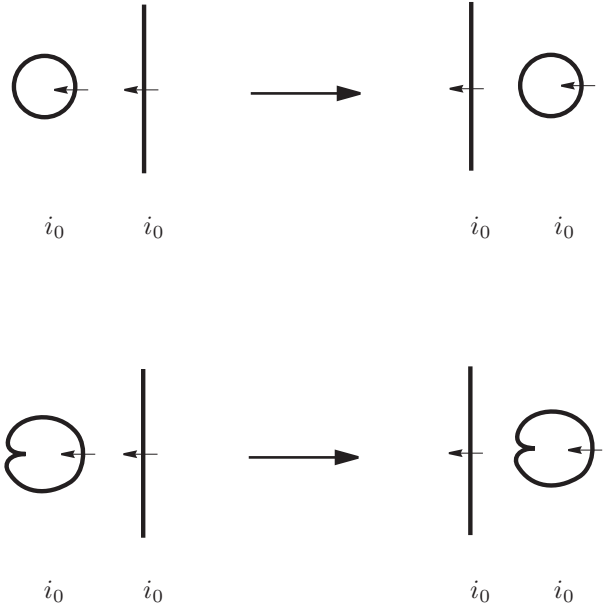}
\caption{Moves replacing a push move}
\label{fig617}
\end{figure}

Then, 
we reach the situation as described
in the leftmost figure of \cite[Fig.~36]{BS1},
where the absolute indices are all equal to
$i_0$, the singular point set image is
topologically embedded, and each component is normally oriented
inward. Note that one of the components
may contain a cusp image.

This completes the proof for the case $n \geq 6$.

When $n=4$, we first follow the procedure as
described in the proof of Lemma~\ref{lem:directed},
taking into account the following points. Note
that in this case $i_0 = 1$.

(1) If we have a cusp, then 
we can use cusp-fold crossing
moves which increase the number of crossings.
Note that they are always-realizable moves.
In our process, these moves are sufficient
for our purpose.

(2) In our process, we create
fold circles of absolute index $i_0+1$
whose images are embedded and are outward-directed.
If we reverse the normal orientations, then
they turn to embedded fold circles of
absolute index $i_0 = 1$, which are inward-directed.
Then, we can safely use the moves as
described in Fig.~\ref{fig617}.

As a result, we get a stable map $f$
with at most one cusp point such that
the fold image of absolute
index $i_0$ is outward-directed 
except for some small embedded inward-directed
circles of absolute index $i_0$.
Then, when we follow the procedure as described
in the proof of Lemma~\ref{lem:elim},
we put the small embedded circles in the region $R$.
Then we can show that the fibers over the points $R$
adjacent to the fold image normally directed 
toward $NP$ are connected.

Consequently, we can follow the same procedure
as in the case $n \geq 6$. This completes the proof.
\end{proof}

\begin{proof}[Proof of Theorem~\textup{\ref{thm1}}]
When $n$ is odd, by applying Lemma~\ref{lem:elim}
successively, we get a desired stable map into $S^2$.

When $n$ is even, by applying Lemma~\ref{lem:elim}
successively, we get a stable map into $S^2$
with at most one cusp such that all the folds
have absolute index $(n-2)/2$. Then, by
applying Lemma~\ref{lem:emb}, we get a desired
stable map.

This completes the proof.
\end{proof}

\begin{rk}\label{rem:Lef}
In the case of dimension $4$ in \cite{BS1},
we have allowed the Lefschetz singularities modeled on
the map germ
$$(z_1, z_2) \mapsto z_1^2 + z_2^2$$
at the origin through $C^\infty$ complex coordinates,
and used always-realizable moves involving
such Lefschetz singularities as well.
In higher even dimensions, say $n = 2m$, we also
have a similar singularity modeled on
$$(z_1, z_2, \ldots, z_m) \mapsto z_1^2 + z_2^2 +
\cdots + z_m^2.$$
It would be an interesting problem to use
this kind of singularities together with folds
and cusps to obtain another proof to the
above results.
\end{rk}

\section{Maps into $\R^2$}\label{section5}

In this section, we prove the following.

\begin{thm}\label{thm2}
Let $M$ be a closed $n$--dimensional manifold with $n \geq 4$, where
$n$ is even.
Then, there exists an image simple
generic map $f : M \to \R^2$ with at most one cusp point.
\end{thm}

\begin{rk}\label{remarkn22}
Theorem~\ref{thm2} holds also for $n=2$; in fact,
for each closed connected surface $M$, we can explicitly
construct a generic map as described in the theorem.
When $M$ is orientable, this is straightforward.
For the non-orientable case, we can first
construct an image simple generic map into $S^2$
by using \cite[Theorem~1.4 (2)]{TY}, for example, 
and then compose it with an appropriate generic map
$S^2 \to \R^2$ to get a desired generic map.
\end{rk}

\begin{proof}[Proof of Theorem~\textup{\ref{thm2}}]
There exists a generic map $f : M \to S^2$ which is image simple
by Theorem~\ref{thm1}.
Let $D_S \cup A \cup D_N = S^2$ be a decomposition
of $S^2$ such that $D_S$ and $D_N$ are disjoint $2$--disks, that $A$
is an annulus, that they are attached along their boundaries,
and that $\Int{D_N} \supset f(S(f))$.
For a regular value $q \in D_S$, set $F = f^{-1}(q)$,
which is an $(n-2)$--dimensional closed manifold.
Note that $f|_{f^{-1}(D_S \cup A)}$ is a trivial smooth
$F$--bundle over $D_S \cup A$ by Ehresmann's fibration theorem
\cite{Eh}. 
In particular, we have
$f^{-1}(A) \cong F \times [-1, 1] \times S^1$.

Let $h : F \to [1, 2]$ be a Morse function such that
the critical values are all distinct. 
Define the Morse function $\varphi: F \times [-1, 1] \to [1, 3]$
by $\varphi(x, t) = h(x) \cos(\pi t/2)+1$, $(x, t) \in 
F \times [-1, 1]$ (refer to \cite[Fig.~40]{BS1}).
Then, we can
construct a stable map $g : M \to \R^2$
in such a way that $g|_{f^{-1}(D_N)} = i_N \circ f|_{f^{-1}(D_N)}$,
$g|_{f^{-1}(D_S)} = i_S \circ f|_{f^{-1}(D_S)}$, and 
$g|_{f^{-1}(A)}$ is the composition of
a trivial $S^1$--family of Morse functions $\varphi \times
\mathrm{id}_{S^1} : (F \times [-1, 1]) \times
S^1 \to [1, 3] \times S^1$
and an appropriate embedding $[1, 3] \times S^1 \to \R^2$,
where $i_N : D_N \to S^2$ and $i_S : D_S \to S^2$ are
appropriate embeddings (see Fig.~\ref{fig612}).
For details, refer to \cite[Proof of Theorem~7.1]{BS1}.

\begin{figure}[htbp]
\centering
\psfrag{F}{$F$}
\psfrag{d1}{$D_N$}
\psfrag{d2}{$D_S$}
\psfrag{fs}{$f(S(f))$}
\includegraphics[width=\linewidth,height=0.6\textheight,
keepaspectratio]{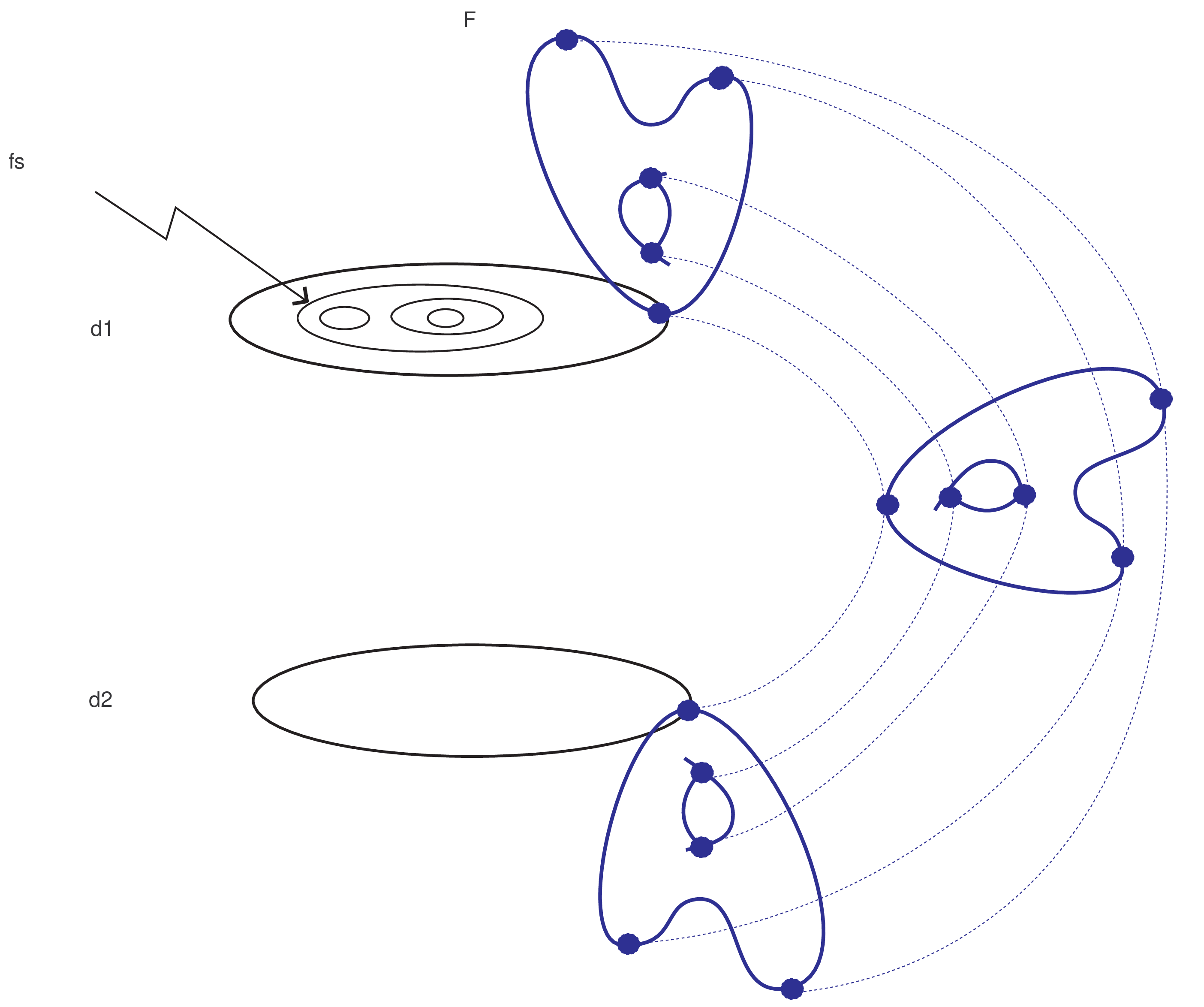}
\caption{Construction of a stable map into the plane}
\label{fig612}
\end{figure} 

The resulting map is an image simple generic map with
at most one cusp. This completes the proof.
\end{proof}

\begin{rk}
The author \cite{Sa3} showed that
a closed orientable $3$--dimensional manifold $M$ admits
an image simple generic map
$f : M \to N$ without cusp points for a surface $N$ if and only if
$M$ is a graph manifold. Recall that a closed orientable
$3$--dimensional manifold is a \emph{graph manifold}
if it is the union of some $S^1$--bundles over compact
surfaces attached along the torus boundaries.
Note that the class of graph manifolds is relatively
small in the set of closed orientable $3$--dimensional
manifolds: for example, hyperbolic $3$--manifolds never
admit such an image simple generic map.

Furthermore, by \cite{SY}, for every ``generic map''
$f : M \to N$ of
a closed oriented $4$--dimensional manifold $M$ into
an orientable $3$--dimensional manifold $N$,
the restriction
$f|_{S(f)}$ has at least $|\sigma(M)|$ 
triple points, where $\sigma(M) \in \Z$ denotes
the signature of $M$.

Furthermore, Gromov \cite{Gromov} showed that for every
$n \geq 3$, there exist a lot of $n$--dimensional
closed manifolds $M$ such that for any
generic map $f : M \to \R^{n-1}$, $f|_{S(f)}$
cannot be a topological embedding.

These results imply that for these dimensions,
for many of the source manifolds,
no map can be modified in such
a way that the restriction to the singular point set
is a topological embedding. 

However, according to our results of this paper,
when $n \geq 4$ is even, a map of a
closed manifold of dimension $n$ into $S^2$ or $\R^2$ can always
be homotopically modified so that it is an embedding
on its singular point set.
\end{rk}

\begin{rk}
In this paper, basically we are studying the
image of the singular point set of a generic map.
For the position of the singular point set
itself in the source manifold, a reader is
referred to the author's previous work \cite{Sa0}
together with \S\ref{section8} of the present paper.
\end{rk}

\section{Relation to open book structures}\label{section6}

In this section, we give an application of the
above results in the odd dimensional case to
open book structures.

\begin{defn}\label{def:openbook1}
Let $M$ be a closed $n$--dimensional manifold.
Suppose that $B$ is a closed submanifold of $M$ of dimension $n-2$ 
with trivial normal bundle, and let $\tau : N(B) \to B \times D^2$
be a trivialization of the normal disk bundle
$N(B)$ of $B$ in $M$, where $N(B)$ is identified with a closed
tubular neighborhood of $B$ in $M$. 
We assume that there exists a smooth
$S^1$--bundle $\varphi : M \setminus B \to S^1$
such that $\varphi|_{N(B) \setminus B}$
coincides with the composition 
\begin{equation}
N(B) \setminus B \spmapright{\tau|_{N(B) \setminus B}} B \times (D^2
\setminus \{0\}) \spmapright{pr_2} D^2 \setminus \{0\} 
\spmapright{r} S^1,
\label{eq:book}
\end{equation}
where $pr_2$ is the projection to the second factor and
$r$ is the radial projection.
Such a structure of $M$ is called an \emph{open book structure}.
Note that for each $\theta \in S^1$,
the closure $F_\theta$ of $\varphi^{-1}(\theta)$ in $M$ coincides
with $\varphi^{-1}(\theta) \cup B$ and is a compact
$(n-1)$--dimensional submanifold of $M$ with boundary $B$.
The submanifold $B$ is called the \emph{binding}
and each $F_\theta$ is called a \emph{page} of the open book
structure.

For some historical account of open book
structures, the reader is referred
Quinn \cite[\S2]{Quinn} and Myers \cite{My}.
\end{defn}

\begin{defn}
Let $M$ be a compact $n$--dimensional manifold with boundary.
Suppose that $B$ is a proper $(n-2)$--dimensional compact
submanifold of $M$ with trivial normal bundle.
Here a submanifold $B$ is \emph{proper}
if $B \cap \partial M = \partial B$ and $B$ intersects
$\partial M$ transversely.
Let $\tau : N(B) \to B \times D^2$
be a trivialization of the normal disk bundle $N(B)$ of $B$
in $M$, where $N(B)$ is identified with a closed
tubular neighborhood of $B$ in $M$ such that $N(B) \cap
\partial M$ can be identified with the tubular neighborhood
$N(\partial B)$ in $\partial M$.
We assume that there exists a smooth
$S^1$--bundle $\varphi : M \setminus B \to S^1$
such that $\varphi|_{N(B) \setminus B}$
coincides with the composition (\ref{eq:book}) as above
and that $\varphi|_{\partial M \setminus \partial B}$
induces an open book structure with binding $\partial B$
on $\partial M$.
Such a structure of $M$ is also called an \emph{open book structure};
$B$ is called the \emph{binding} and the closure of each
fiber of $\varphi$ is called a \emph{page}.
Note that the restriction of the structure to
$\partial M$ gives rise to an open book structure
of $\partial M$ in the sense of Definition~\ref{def:openbook1}. 
In this case, we say that
the open book structure on the boundary \emph{extends}
through $M$.
\end{defn}

As an application of our results, we give
a new proof to the following result originally due
to Alexander \cite{Alex} for $n=3$
and to Quinn \cite{Quinn} for $n \geq 5$
(see also \cite{Lawson, My, Tamura, Wink}).

\begin{thm}\label{thm:ob}
Let $M$ be a compact connected $n$--dimensional manifold with $n \geq 3$ odd
such that an open book structure on $\partial M$ is given.
Then, it always extends through $M$.
\end{thm}

\begin{proof}
First, let us construct a certain generic map $f : M \to \R^2$
compatible with the open book structure on $\partial M$
such that $f(M)$ coincides with the unit disk $\Delta
\subset \R^2$, as follows.

Let $bB \subset \partial M$ be the binding,
$b\tau : N(bB) \to bB \times D^2$ the trivialization
of the normal disk bundle of $bB$ in $\partial M$, and
$b\varphi : \partial M \setminus bB \to S^1$ 
the fibration associated with the open book structure on $\partial M$.
Let $f|_{N(bB)} : N(bB) \to \R^2$ be defined by the
composition 
$$N(bB) \spmapright{b\tau} bB \times D^2 \spmapright{pr_2}
D^2 \spmapright{\iota} \R^2,$$
where $pr_2$ is the projection to the second factor
and $\iota : D^2 \to \R^2$ is an embedding onto $\Delta$.
We further define $f|_{\partial M \setminus \Int{N(bB)}} : 
\partial M \setminus \Int{N(bB)} \to \partial \Delta$
as the composition of $b\varphi|_{\partial M \setminus \Int{N(bB)}}
: \partial M \setminus \Int{N(bB)} \to S^1$ and an
appropriate diffeomorphism $S^1 \to \partial \Delta$ in such a way
that $f|_{\partial M}$ thus constructed gives a smooth map
$\partial M \to \Delta \subset \R^2$.
Here, strictly speaking, $M$ is considered
to be a manifold with corners along $\partial N(bB)$.
Note that this map has no singular points and that
it can naturally be extended to a submersion, denoted by
$f$, near $\partial M$ so that $f^{-1}(\partial \Delta)
= \partial M \setminus \Int{N(bB)}$.

As $\Delta$ is contractible, we can further
extend the map $f$ over the entire manifold $M$
as a smooth map and then approximate it by a generic map, still
denoted by $f$. We may assume that $f$
satisfies the conditions (1) and (2) of Definition~\ref{def:stable}.
As it is a submersion near $\partial M$,
we see that $S(f)$ is contained in $\Int{M}$
and we may further assume that its image $f(S(f))$ sits
inside $\Int{\Delta}$ and that $f(M) = \Delta$.

In the following, we will deform $f$ successively
for our purpose and will continue
to denote the resulting maps by the same symbol $f$.

We see that the number of cusps of $f$ is even,
since the Euler characteristic of the pair
$(M, F)$ vanishes, where $F$ is a fiber of $f$
over a point in $\partial \Delta$.
Then, by the same procedure as in Levine's cusp
elimination theorem \cite{Lev}, we can eliminate the
cusps of $f$ by pairs, since $M$ is connected. 
So, we may assume that
$f$ has no cusps.

Now, we apply Lemma~\ref{lem:directed}, taking
into account the comments in
Remark~\ref{rk:directed} (1) and (2),
first to the definite fold image.
As a result, we get a generic map
$M \to \R^2$, denoted by $f$ again,
such that $f(M)$ coincides with $\Delta$,
the definite fold image is outward directed,
and that the other fold images lie
inside the interior of a $2$--disk, say $\Delta_1$,
contained in the innermost region.
Then, we apply Lemma~\ref{lem:directed}
to $f$ restricted to $f^{-1}(\Delta_1)$
to get a new generic map $f$ such that the fold
image of absolute index $1$ is contained
in $\Int \Delta_1$, that it is outward directed,
and that the other fold images lie
inside the interior of a $2$--disk, say $\Delta_2$,
contained in the innermost region.
Repeating this procedure, we get a generic
map $f : M \to \R^2$ such that $f(M) = \Delta$,
that the fold images of absolute indices $i$ with
$0 \leq i < \lfloor (n-1)/2 \rfloor$ are all directed outward,
that their positions are consistent with the order of absolute indices
(smaller indices in the outer position),
and that the fold image of absolute index $i_0 = \lfloor (n-1)/2 \rfloor$
lies in the interior of a $2$--disk $\Delta_{i_0}$
contained in the innermost region.
At this stage, we should note that the fold image of absolute
index $i_0$ may not be directed outward.
We may continue to assume that $f$ satisfies
the conditions (1) and (2) of Definition~\ref{def:stable}.

Note that as $n$ is odd, even if we reverse the normal
orientation of the fold image of absolute index
$i_0$, the index does not change.

Now, let us first consider the case $n \geq 5$.
We apply the always-realizable moves as depicted in
Figs.~\ref{fig619}--\ref{fig622}.
We may assume that the disks $\Delta_0, \Delta_1,
\ldots, \Delta_{i_0}$ are concentric, i.e., their
centers all coincide with the origin $\mathbf{0}
\in \R^2$.
We denote by $\pi : \Delta \setminus \{\mathbf{0}\} \to S^1$
the radial projection.
Let $A$ be the annulus $\Delta_{i_0} \setminus \Int{\Delta_{i_0+1}}$,
where $\Delta_{i_0+1}$ is a concentric disk in $\Int{\Delta_{i_0}}$.
We may assume that the fold image of absolute index $i_0$
lies in $\Int{A}$ and that $\pi \circ f|_{S_{i_0}} :  S_{i_0} \to S^1$
is an $S^1$--valued  Morse function with distinct
critical values, where $S_{i_0}$ is the set of folds
of absolute index $i_0$ of $f$.
Then, for each of the critical points, we apply
the type II moves as depicted in Fig.~\ref{fig619}
so that the images of the critical points by $f$
lie inside $\Int{\Delta_{i_0+1}}$.
(Before doing so, we arrange $f(S_{i_0})$
so that all the crossings are pushed into the
regions in $A$ outside of those regions in which
$f(S_{i_0})$ is depicted
in the figure.)
Note that in our case, only the top type II move
(from right to left)
in Fig.~\ref{fig605} is allowed and that
the bottom type II move is not. Note also that
as $n$ is odd,
even if we reverse the normal orientation for
a fold image of absolute index $i_0 = (n-1)/2$,
the absolute index does not change.

\begin{figure}[htbp]
\centering
\psfrag{A}{$A$}
\psfrag{D}{$\Delta_{i_0+1}$}
\psfrag{D0}{$\Delta_{i_0}$}
\psfrag{o}{$\textbf{0}$}
\psfrag{S}{$f(S_{i_0})$}
\includegraphics[width=\linewidth,height=0.5\textheight,
keepaspectratio]{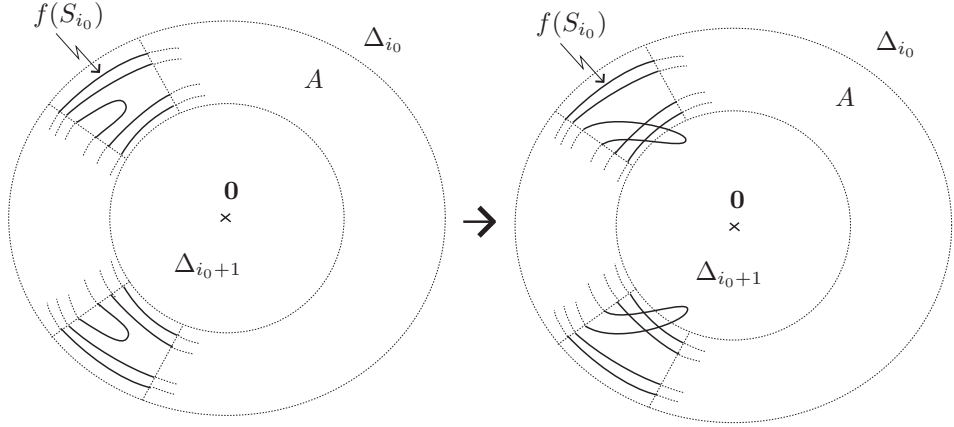}
\caption{Pushing the local minima and maxima into $\Int{\Delta_{i_0+1}}$.}
\label{fig619}
\end{figure}

We orient $S^1 = \partial \Delta$ counter clockwise, for the
moment, then with respect to this orientation, the numbers of
local minima and maxima of the $S^1$--valued Morse function
$\pi \circ f|_{S_{i_0}}$ are the same. Then, by further applying the
type II moves, we may assume that the images
of the critical points by $f$ lie in pairs, ``face-to-face''
as in the first figures of Figs.~\ref{fig620} and \ref{fig621}.
Here, we can show that the fibers of $f$ over the points
in $\Delta_{i_0+1}$ are connected by the same argument
as in the proof of Lemma~\ref{lem:elim}, since $n \geq 5$.
Then, we apply always-realizable moves as
depicted in Fig.~\ref{fig620} for each pair
of the images of a local minimum and a local maximum.
More precisely, we first apply type I moves to create
4 cusps. Then, we apply the cusp merge moves to the 2 pairs
of cusps, which are allowed, since the fibers are connected.
If we reverse the normal orientation of the newly born
index--$(i_0+1)$ circle, then it turns into a fold image circle
of absolute index $i_0-1$. Then, by applying
the moves as depicted in Fig.~\ref{fig611} with
$i_0$ replaced by $i_0-1$, we can turn this circle
into a circle of absolute index $i_0$. (Note that
$i_0-1 \leq (n-3)/2$, which guarantees these moves.)
Then, by applying the type II move once, we reach
the final configuration in Fig.~\ref{fig620}.

\begin{figure}[t]
\centering
\psfrag{i}{$i_0$}
\psfrag{i1}{$i_0+1$}
\psfrag{i2}{$i_0-1$}
\includegraphics[width=\linewidth,height=0.5\textheight,
keepaspectratio]{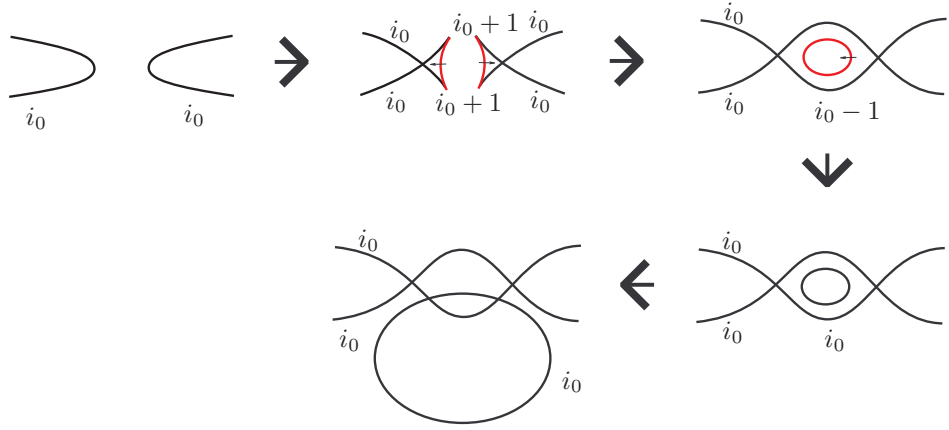}
\caption{A series of moves near a pair of critical point images}
\label{fig620}
\end{figure}

By applying this procedure for each pair of critical point images,
we reach the configuration as depicted in the figure on the
right hand side in Fig.~\ref{fig621}.

\begin{figure}[t]
\centering
\psfrag{o}{$\mathbf{0}$}
\psfrag{D}{$\Delta_{i_0+1}$}
\includegraphics[width=\linewidth,height=0.5\textheight,
keepaspectratio]{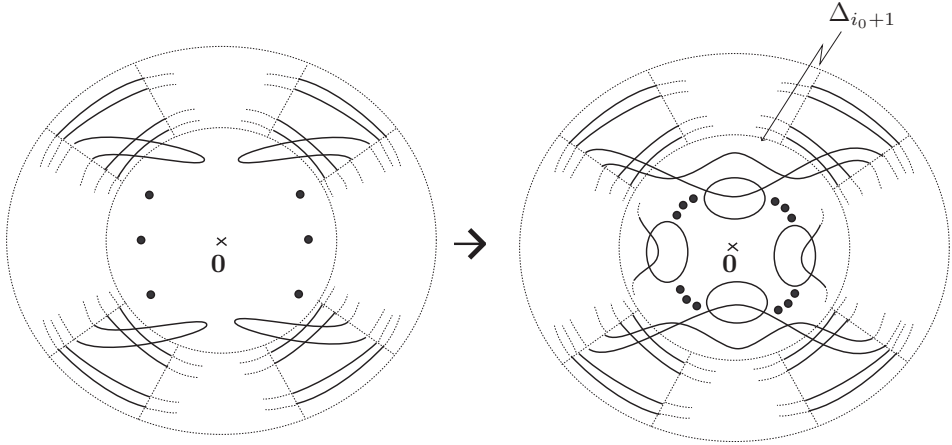}
\caption{Applying the series of moves as in Fig.~\ref{fig620}
for each pair}
\label{fig621}
\end{figure}

Then, inside $\Int{\Delta_{i_0+1}}$, 
we move the circles using the type II moves
several times as depicted in Fig.~\ref{fig622}.
We perform the moves in such a way that for
the resulting map, still denoted by $f$,
$\pi \circ f|_{S_{i_0}} : S_{i_0} \to \partial \Delta = S^1$
is a submersion.

\begin{figure}[t]
\centering
\psfrag{o}{$\mathbf{0}$}
\includegraphics[width=\linewidth,height=0.5\textheight,
keepaspectratio]{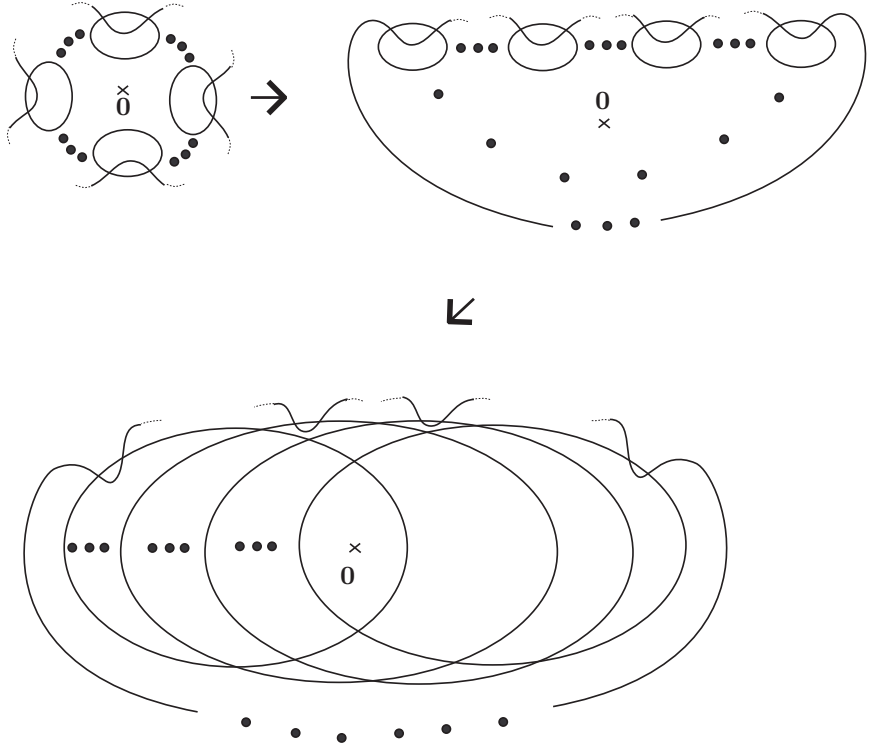}
\caption{Arranging the circles toward the center $\mathbf{0}$}
\label{fig622}
\end{figure}

Let $\Delta_{i_0+2}$
be a small closed $2$--disk centered at $\mathbf{0}$
disjoint from
$f(S(f))$. Then, as $f$ restricted to $f^{-1}(\Delta_{i_0+2})$
and that to $f^{-1}(\Delta_{i_0+2}) \cap \partial M$
are submersions, we see that $f|_{f^{-1}(\Delta_{i_0+2})}
: f^{-1}(\Delta_{i_0+2}) \to \Delta_{i_0+2}$
is a smooth fiber bundle, and it is trivial, since $\Delta_{i_0+2}$
is contractible. Furthermore, setting $B = f^{-1}(\mathbf{0})$,
we see that $\partial B = bB$ and that the composition
$$M \setminus B \spmapright{f} \Delta \setminus \{\mathbf{0}\} \spmapright{\pi} S^1$$
is a smooth fiber bundle, since $\pi \circ f$ restricted to $S(f)$
is a submersion.
Hence, this gives an open book structure for $M$.

Precisely speaking, this open book structure
does not induce the given open book structure
on the boundary $\partial M$. However,
making $N(bB)$ smaller, we see that the
open book structure induced
on $\partial M$ is
naturally isomorphic to the given one.
This means that the given open book structure
on $\partial M$ extends through $M$.
This completes the proof for the case $n \geq 5$.

When $n=3$, we basically follow the same
procedure.
However, the fibers of $f$ over the points
in $\Delta_{i_0+1}$ may not be connected,
even if the $3$--manifold $M$ is connected.
Therefore, the cusp merges as in Fig.~\ref{fig620}
may not be applicable. 
In this case, we can find a pair of parallel embedded curves
$\gamma_1$ and $\gamma_2$
inside $M$ connecting the relevant cusp points
such that their interiors avoid $S(f)$ and that
their images by $f$ are immersed.
Then, by imitating Alexander's argument \cite{Alex},
we may further arrange $\gamma_1$ and $\gamma_2$ so that their $f$--images
wind around the origin in the same direction.
This is possible, since the folds have absolute indices
$0$ or $1$ and those of absolute index $0$ are outward
directed.
We can, then, merge the $2$ pairs of cusps along
$\gamma_1$ and $\gamma_2$, using fold-cusp crossings as well.
See Fig.~\ref{fig625}.

\begin{figure}[t]
\centering
\psfrag{o}{$\mathbf{0}$}
\psfrag{D}{$\Delta_{i_0+1}$}
\includegraphics[width=\linewidth,height=0.8\textheight,
keepaspectratio]{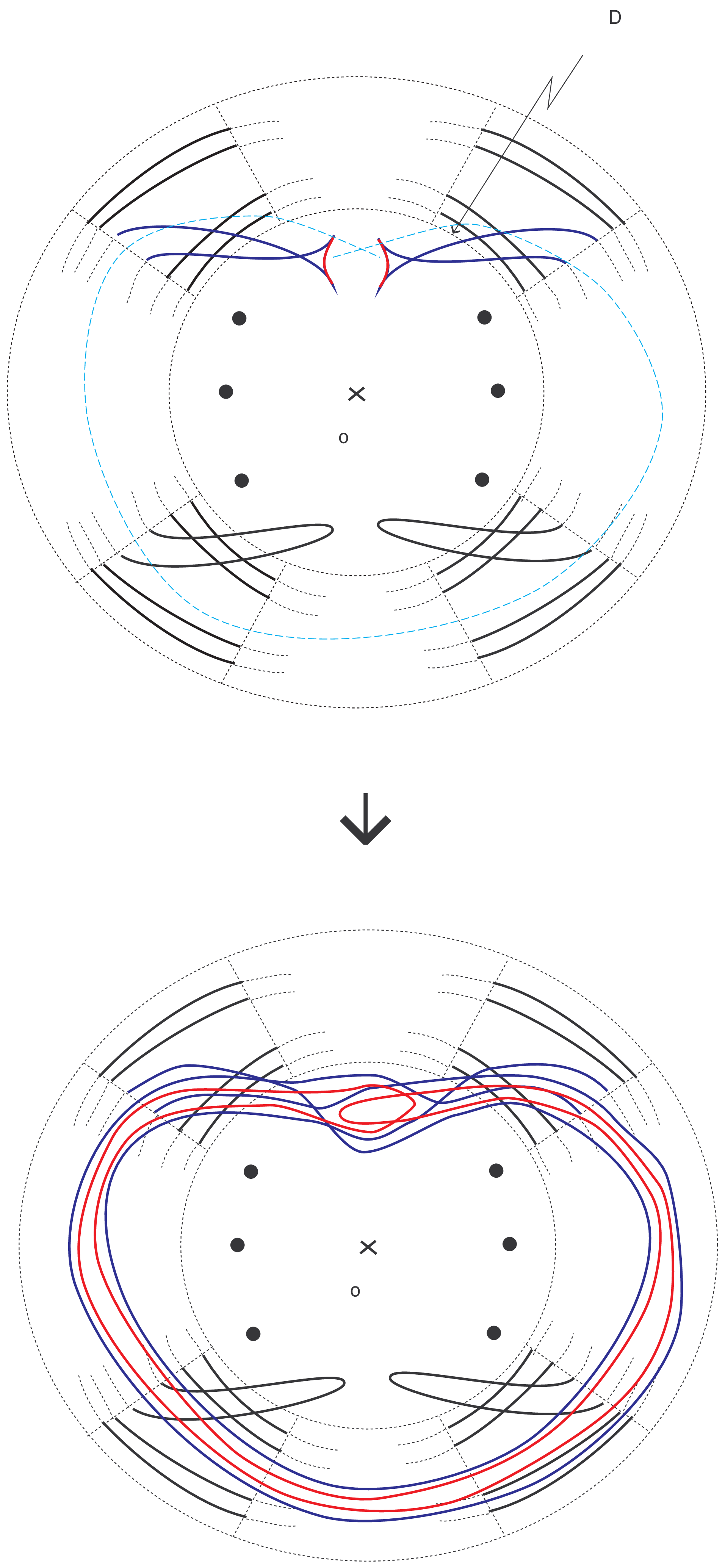}
\caption{Cusp merges along curves in $M$}
\label{fig625}
\end{figure}

By choosing the parallel curves $\gamma_1$ and $\gamma_2$ 
appropriately, 
we may assume that the newly created circle
of definite fold image bounds an immersed disk.
We can, then, shrink the immersed disk
so that we finally have an embedded one, by
sliding the definite fold image circle
toward the small disk embedded inside
the immersed disk. Then we can
safely apply the move as in Fig.~\ref{fig611}
in order to replace the inward-directed
definite fold image circle
with that of absolute index $1$.
Then, we can pull it toward the origin
as in the last step of Fig.~\ref{fig620}.
The only difference is that now the circle may
intersect with several parallel fold images
of absolute index $1$, which causes no problem.

We perform these procedures for each facing pair
of the images of a local maximum and a local minimum.
Finally, we perform the moves as in Fig.~\ref{fig622}.
Then we get to the same conclusion as in the case of
$n \geq 5$. This completes the proof.
\end{proof}

\begin{rk}
(1) Note that Theorem~\ref{thm:ob} also covers the case
with $\partial M = \emptyset$: just ignore
the boundary open book structures in the proof.

(2) In the case $n=3$, we do not know if
we can make the binding connected by modifying our proof.
\end{rk}

\section{Elimination of indefinite folds of certain indices}\label{section7}

In this section, using the constructions of the previous
section for the odd dimensional case, we show that
if $M$ is an odd dimensional closed manifold
which is $k$--connected for a certain positive integer $k$,
then we can eliminate folds of absolute indices
$i$ with $1 \leq i \leq k$ for maps into $\R^2$.

Let us start with the following.

\begin{lem}\label{lemma:p1}
For the extended open book structure on $M$ constructed
in the above proof of Theorem~\textup{\ref{thm:ob}}, the binding $B$
is connected, and the inclusion $B \to M$ induces an
isomorphism $\pi_1(B) \cong \pi_1(M)$ between the
fundamental groups if $n \geq 7$.
\end{lem}

\begin{proof}
The binding $B$ is connected as shown in the proof
of Theorem~\textup{\ref{thm:ob}}.
Let us consider the fundamental group.

Take the segment $J = [0, 1] \times \{0\} \subset \Delta$.
We see that $F = f^{-1}(J)$ coincides with a page
of the open book structure of $M$ and that
$f|_F : F \to J = [0, 1]$ is a Morse function
such that $(f|_F)^{-1}(0) = B = \partial F$.
Furthermore, the indices of the critical points
are all greater than or equal to $i_0 = (n-1)/2 \geq 3$.
Hence, the inclusion map $B \to F$ induces an isomorphism
between the fundamental groups.

Let us decompose $\Delta$, the unit disk in $\R^2$,
into two half disks as $\Delta = \Delta_+ \cup \Delta_-$, where
$\Delta_+ = \{(x, y) \in \Delta\,|\, y \geq 0\}$ and
$\Delta_- = \{(x, y) \in \Delta\,|\, y \leq 0\}$.
Then, we see that the inclusions
$F \to f^{-1}(\Delta_{\pm})$ are homotopy equivalences.
Now, the intersection $X = f^{-1}(\Delta_+) \cap f^{-1}(\Delta_-)$
is diffeomorphic to the union of $F$ and its copy
attached along $B$ and hence we see that the inclusion
$B \to X$ induces an isomorphism between the
fundamental groups by the van Kampen theorem.
Finally, by the van Kampen theorem applied to
$M = f^{-1}(\Delta_+) \cup_X f^{-1}(\Delta_-)$,
we see that the inclusion $B \to M$
induces an isomorphism between the fundamental groups.
This completes the proof.
\end{proof}

\begin{thm}\label{thm:ob2}
Let $M$ be a closed simply connected
$n$--dimensional manifold with $n \geq 7$ odd.
Then, there exists a stable map $f : M \to \R^2$
without cusp points such that the absolute
indices $i$ of the folds satisfy either $i = 0$
or $2 \leq i \leq (n-1)/2$.
\end{thm}

\begin{proof}
Let us consider the generic map $f : M \to \R^2$
constructed in the proof of Theorem~\ref{thm:ob}.
This is a generic map without cusp points; however, it
may contain fold points of absolute index $1$.
We will eliminate such fold points by homotopy as follows.

As in the proof of Lemma~\ref{lemma:p1},
we consider the Morse function $g = f|_F : F \to
[0, 1]$, where $F$ is a compact $(n-1)$--dimensional
manifold with boundary $\partial F = B$.
By the construction of $f$, we see that
the indices $i$ of the critical points satisfy
$(n-1)/2 \leq i \leq n-1$. On the other hand,
by our assumption and Lemma~\ref{lemma:p1}, 
$B$ and $F$ are simply connected.

If $g$ has two or more critical points of
index $n-1$, i.e.\ local maxima,
then there exists a critical point of index $n-2$
in such a way that its dual handle connects
two dual $0$--handles, since $F$ is connected. 
Therefore, either of the two
critical points of index $n-1$ and the index $n-2$
critical point make a canceling pair.
(Note that by reversing the orientation of $J = [0, 1]$,
they make a canceling $0$--$1$ pair of critical points.)
Then by \cite[Part~I, Chapter~V, Proposition~1.1]{HW},
we can cancel them leaving a pair of cusps
of absolute index $0$. By repeating this procedure,
we may assume that $g$ has only one maximal point
(refer to Fig.~\ref{fig623}).

\begin{figure}[htbp]
\centering
\psfrag{0}{$0$}
\psfrag{1}{$1$}
\psfrag{J}{$J$}
\includegraphics[width=\linewidth,height=0.7\textheight,
keepaspectratio]{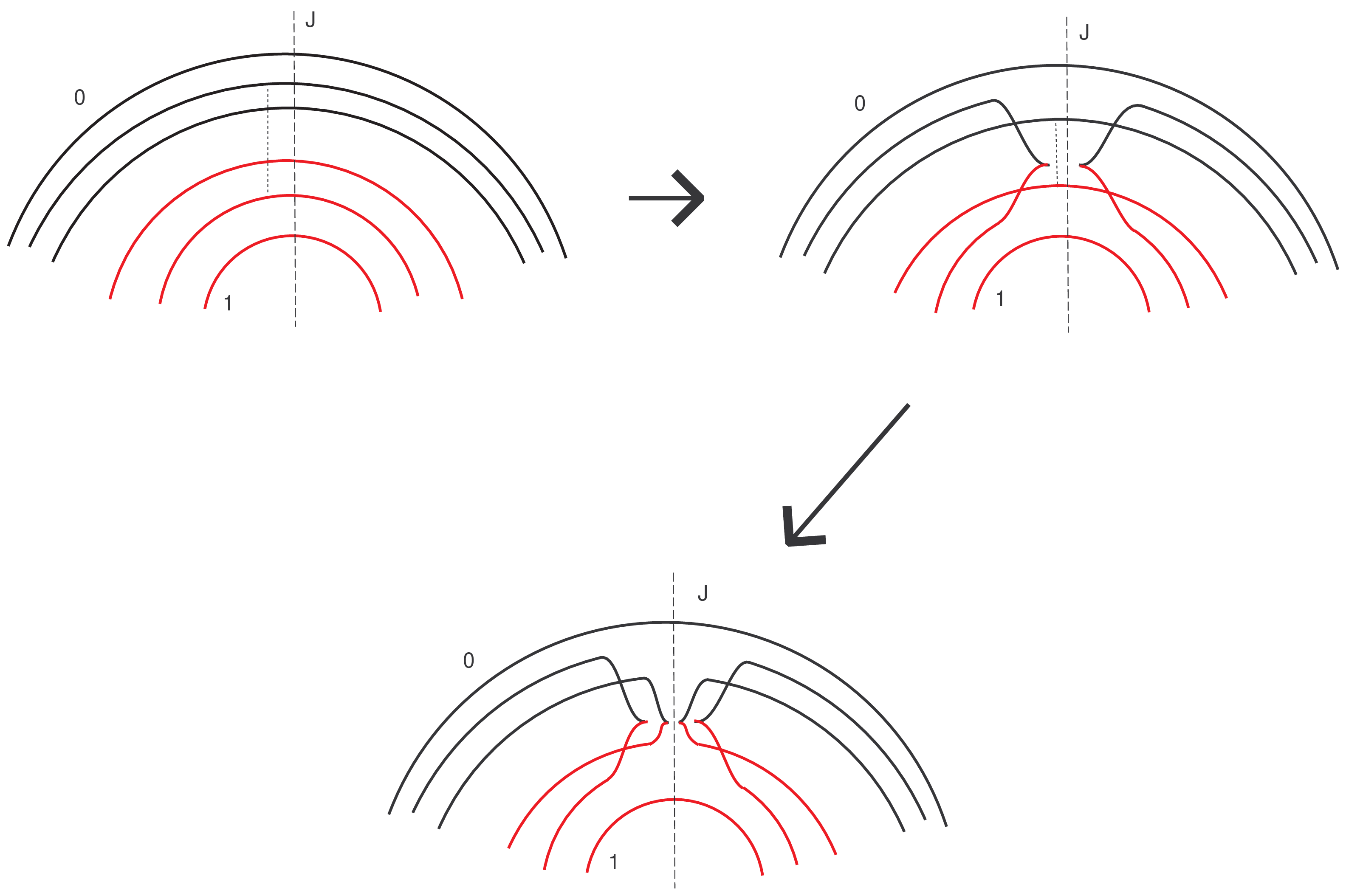}
\caption{Canceling $0$--$1$ pairs of critical points}
\label{fig623}
\end{figure}

Let $U$ be a small open regular neighborhood of $J$ in $\R^2$.
Note that $\Delta \setminus U$ is naturally identified with
$[0, 1] \times [-1, 1]$, where the first factor corresponds to
a radial outward arc and the second factor
corresponds to the direction of
$\partial \Delta = S^1$. 
Then, $f$ restricted to $f^{-1}(\Delta \setminus U)$ 
is considered to be a smooth $1$--parameter family of
functions from $F$ to $[0, 1]$ parametrized by $[-1, 1]$.
Now, by applying an argument similar to that
in the proof of \cite[Part~I, Chapter~V, Proposition~3.5]{HW},
we can homotopically modify $f$ so that
it has only one fold circle of absolute index $0$
and that it is embedded into $\R^2$ in the
outermost place with the expense of creating
some fold circles of absolute index $2$.

Since $F$ is simply connected of dimension
greater than or equal to $6$, 
by \cite[Theorem~8.1]{M1}, we can cancel
each critical point of index $1$ with the expense
of creating critical points of index $3$.
Then, again by the same argument as 
in the proof of \cite[Part~I, Chapter~V, Proposition~3.5]{HW},
we can homotopically modify $f$ so that
it has no fold points of absolute index $1$ and that
it has only one fold circle of absolute index $0$
embedded into $\R^2$ in the
outermost place.
This completes the proof.
\end{proof}

Now, as a small digression,
let us consider the topology of manifolds
admitting a generic map into $\R^2$ without fold
of absolute index $1$.

\begin{lem}\label{lem:free}
Let $M$ be a closed connected $n$--dimensional manifold,
$n \geq 3$, and $f : M \to \R^2$ a generic map without fold
of absolute index $1$. Then, the fundamental group
$\pi_1(M)$ of $M$ is a free group of finite rank.
\end{lem}

\begin{proof}
Let us consider the quotient space $W_f$ as follows.
For two points $x, x' \in M$, we define $x \sim x'$ if
$f(x) = f(x')$ and $x, x'$ lie on the same
connected component of $f^{-1}(f(x)) = f^{-1}(f(x'))$.
Then $W_f$ is defined to be the quotient space $M/\!\sim$.
Note that for the quotient map $q_f : M \to W_f$,
there exists a unique continuous map $\bar{f} : W_f \to \R^2$
such that $f = \bar{f} \circ q_f$.

We may assume that $f$ is a stable map. Then, 
as $f$ does not have folds of absolute index $1$,
we can show that $W_f$ has the structure of a compact
connected smooth $2$--dimensional manifold
(for example, see \cite[Theorem~9.1.7]{Wr} or \cite{Sa4}).
Note that $q_f$ restricted to the definite fold locus
is a diffeomorphism onto
the boundary $\partial W_f$.
As we can see easily that $\bar{f} : W_f \to \R^2$ is an immersion,
$W_f$ is orientable and has nonempty
boundary. Let $\alpha_1, \alpha_2, \ldots, \alpha_r$
be disjoint properly embedded arcs in $W_f$
such that $Y = W_f \setminus \left(\cup_{s=1}^r \Int{N(\alpha_s)}\right)$
is diffeomorphic to the $2$--dimensional disk (with $4r$ corner
points), where $N(\alpha_s)$ are disjoint small regular neighborhoods
of $\alpha_s$ in $W_f$. We may assume that
the images of the cusp points by $q_f$ and
the crossings of $q_f|_{S(f)}$ lie in $\Int{Y}$,
that $q_f|_{S(f)}$
is transverse to each $\alpha_s$, and 
that each connected component of $q_f(S(f)) \cap N(\alpha_s)$
is a properly embedded arc in $N(\alpha_s)$ intersecting with
$\alpha_s$ at exactly one point (see the top figure
of Fig.~\ref{fig624}).

\begin{figure}[htbp]
\centering
\psfrag{N1}{$N(\alpha_1)$}
\psfrag{N2}{$N(\alpha_2)$}
\psfrag{Nr}{$N(\alpha_r)$}
\psfrag{b}{$\beta$}
\psfrag{q}{$q_f(S(f))$}
\psfrag{Y}{$Y$}
\psfrag{W}{$W_f$}
\psfrag{h}{$h$}
\psfrag{R}{$\R$}
\includegraphics[width=\linewidth,height=0.8\textheight,
keepaspectratio]{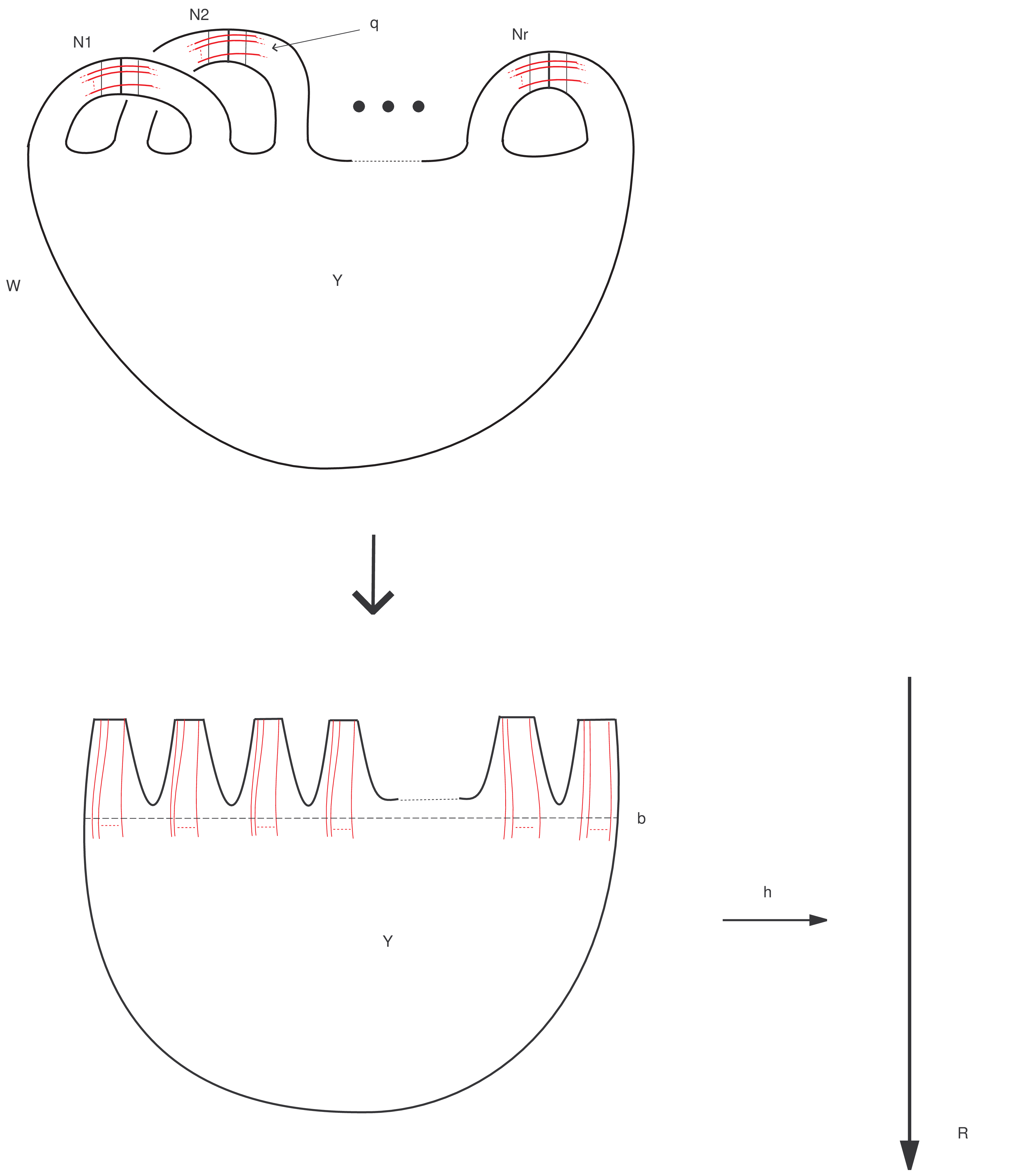}
\caption{Cutting open the quotient space}
\label{fig624}
\end{figure}

Let us cut open $W_f$ along $N(\alpha_s)$, $s = 1, 2,
\ldots, r$, so that we obtain the bottom figure of
Fig.~\ref{fig624}. Let $\beta$ be a properly
embedded arc in $Y$ as in the figure.
Let us consider the structure of $q_f^{-1}(\beta)$.
First, note that for each $s$, $M_s = q_f^{-1}(\alpha_s)$
is a closed connected $(n-1)$--dimensional manifold and that
$q_f|_{M_s} : M_s \to \alpha_s \cong [-1, 1]$
is a Morse function with exactly one minimum and one maximum
and without critical points of index $1$ or $n-2$.
Hence $M_s$ is simply connected.
Let $h : Y \to \R$ be a smooth function which
corresponds to the negative of the height function
as depicted in the bottom figure of Fig.~\ref{fig624}.
Then we may assume that $h \circ q_f : q_f^{-1}(Y) \to \R$
is a Morse function. Note that its critical points
correspond to those of $h \circ q_f$ restricted to
the singular point set $S(f) \cap q_f^{-1}(Y)$,
and we may assume that the cusp points are not critical points.
Since $f$ does not have folds of absolute index $1$,
we see that the critical points on indefinite folds
have indices different from $0, 1, n-1$ or $n$.
Furthermore, the function
has a unique critical point of index $n$ together with
exactly $r-1$ critical
points of index $1$, each of which corresponds to
a $1$--handle connecting $M_{s_1}$ and $M_{s_2}$, $s_1 \neq s_2$.
Hence, $q_f^{-1}(\beta)$ is diffeomorphic to the
connected sum $M_1 \sharp M_2 \sharp \cdots
\sharp M_r$, which is simply connected.

By considering the Morse function $h \circ q_f : q_f^{-1}(Y) \to \R$,
we see that $q_f^{-1}(Y)$
is also simply connected. Thus, the original manifold $M$
is obtained from the simply connected manifold
$q_f^{-1}(Y)$ by attaching simply connected manifolds
$q_f^{-1}(N(\alpha_s)) \cong M_s \times [-1, 1]$
along their boundaries, $s = 1, 2, \ldots, r$.
Hence, by an argument using the van Kampen theorem,
we see that the fundamental group of $M$
is a free group of rank $r$.
This completes the proof.
\end{proof}

\begin{rk}
The above lemma implies that if we drop the
condition that $M$ should be simply connected
in Theorem~\ref{thm:ob2}, the conclusion
does not hold in general: if $\pi_1(M)$
is not a free group, a generic map as described
in the theorem
does not exist. On the other hand,
if $\pi_1(M)$ is a nontrivial free
group, then we do not know if
a generic map as described
in the theorem exists or not.
\end{rk}

Let us now consider elimination of folds with
higher absolute indices. For this, we need the
following.

\begin{lem}\label{lemma:p2}
For the extended open book structure on $M$ constructed
in the proof of Theorem~\textup{\ref{thm:ob}}, 
if $M$ is $k$--connected for an integer $k$ with
$1 \leq k \leq (n-5)/2$, then
the binding $B$ and the page $F$ are also $k$--connected.
\end{lem}

\begin{proof}
In the following, we use the same
notation as in the proof of Lemma~\ref{lemma:p1}.

First note that $F$ is obtained by
attaching handles of index $\geq (n-1)/2$ to
$B \times [0, 1]$, and hence that the
inclusion map $B \to F$ induces isomorphisms
in homology of dimensions $< (n-3)/2$.

Set $X_{\pm} = f^{-1}(\Delta_{\pm})$ and
$X = f^{-1}(\Delta_+) \cap f^{-1}(\Delta_-)$.
Then,
we have the following Meyer--Vietoris exact sequence
of homology with coefficients in $\Z$:
\begin{equation}
\cdots \to H_{i+1}(M) \to H_i(X) \to H_i(X_+) \oplus
H_i(X_-) \to H_i(M) \to \cdots.
\label{exact1}
\end{equation}
We also have the exact sequence
$$\cdots \to H_{i+1}(X) \to H_i(B) \to H_i(F)
\oplus H_i(F) \to H_i(X) \to \cdots.$$
As $H_i(B) \to H_i(F)
\oplus H_i(F)$ is injective for $i < (n-3)/2$, we
have the short exact sequence
$$0 \to H_i(B) \to H_i(F)
\oplus H_i(F) \to H_i(X) \to 0.$$
This implies that the inclusion
induces an isomorphism $H_i(F) \to H_i(X)$
for $i < (n-3)/2$.
Thus, the inclusions $X \to X_{\pm} \simeq F$ also
induce isomorphisms in homology
for dimensions $< (n-3)/2$.
Hence, by (\ref{exact1}) we have the short exact sequence
$$0 \to H_i(X) \to H_i(X_+) \oplus
H_i(X_-) \to H_i(M) \to 0$$
for $i < (n-3)/2$.
This implies that the inclusion
induces an isomorphism $H_i(F) \to H_i(M)$
for $i < (n-3)/2$.
Since $B$ and $F$ are simply connected
by Lemma~\ref{lemma:p1} and $M$ is $k$--connected
with $k < (n-3)/2$, we see that
$B$ and $F$ are also $k$--connected
by the Hurewicz theorem.
This completes the proof.
\end{proof}

Then, using the above lemma and an argument
similar to that in the proof of Theorem~\ref{thm:ob2},
we get the following.

\begin{thm}\label{thm:ob3}
Let $M$ be a closed simply connected
$n$--dimensional manifold with $n \geq 7$ odd.
Suppose $M$ is $k$--connected
for an integer $k$ with
$1 \leq k \leq (n-5)/2$.
Then, there exists a stable map $f : M \to \R^2$
without cusp points such that 
$f$ has no fold points of absolute index $i$ with
$1 \leq i \leq k$.
\end{thm}

\begin{proof}
By Theorem~\ref{thm:ob2}, we can eliminate
folds of absolute index $1$ with the expense
of creating folds of absolute indices $2$ and $3$.
As $F$ is $k$--connected, we can use the same
argument to successively eliminate folds
of absolute indices $2, 3, \ldots, k$ with
the expense of creating folds of absolute
indices $k+1$ and $k+2 \leq (n-1)/2$.
This completes the proof.
\end{proof}

\begin{rk}
In the proofs of Theorems~\ref{thm:ob2}
and \ref{thm:ob3}, we have constructed
generic smooth maps of $M$ into $\R^2$
without cusps which are compatible with
the open book structure constructed
in the proof of Theorem~\ref{thm:ob}.
\end{rk}

\begin{rk}
We do not know if the condition $1 \leq k \leq (n-5)/2$
can be weakened or not. For example,
we do not know if a $1$--connected
$5$--dimensional manifold always admits a generic map into $\R^2$
without cusps and with only folds of absolute
indices $0$ or $2$.
\end{rk}

\begin{rk}
In Theorem~\ref{thm:ob3}, the condition that
$M$ should be $k$--connected is not a necessary
condition for the existence of
a stable map $f : M \to \R^2$
without cusp points such that 
$f$ has no fold points of absolute index $i$ with
$1 \leq i \leq k$.
As the proof of Lemma~\ref{lem:free} shows,
even if such an $f$ exists, $M$ may have
nontrivial free fundamental group.
For example, there exist a plenty of explicit
examples of manifolds with nontrivial
free fundamental groups which admit generic maps
into $\R^2$ with only folds of absolute index $0$
(see \cite{Sa93}, for example).
\end{rk}

\begin{rk}
We do not know if in Theorem~\ref{thm:ob3}
the condition that $n \geq 7$ is odd and
$1 \leq k \leq (n-5)/2$ is essential or not.
\end{rk}

\section{Singular point set}\label{section8}

In this section, we consider the topological
position of the singular point set of a generic
map, instead of its image. Such study was
done in the author's previous paper \cite{Sa0};
however, in some cases, the result was not
correctly stated. In this section, we give
a correction to the result.

More precisely, in \cite{Sa0}, the author has studied the position
of the singular point set $S(f)$ of a generic
map $f : M \to N$ of a closed $n$--dimensional
manifold, $n \geq 3$, into a surface $N$.
However, when $n$ is even, the main result is
not correctly stated and the proof needs also to be modified
appropriately, as follows.

Let $f : M \to N$ be a generic map of a closed
$n$--dimensional manifold into an oriented surface $N$,
where we assume that $n \geq 2$ is even.
In this case using the map $f|_{S(f)}$, we can
naturally orient $S(f)$ as follows.
Let $C$ be a component of $S(f) \setminus \{\mbox{cusps}\}$
with absolute index $i$. Recall that $f|_C : C \to N$
is a normally oriented immersion. For a point $p \in C$,
we orient $T_pC$ so that the ordered basis of $T_{f(p)}N$
consisting of a normal orientation and the image
by $df_p$ of a nonzero vector in $T_pC$
consistent with the orientation gives the orientation
(or opposite orientation) of $T_{f(p)}N$ when $i$ is even
(resp.\ odd). Then, with the help of Fig.~\ref{fig601},
we see that this gives a unique orientation
for each component of $S(f)$.
In the rest of this section, we always consider
this orientation for $S(f)$.

Note that such an orientation cannot be defined
canonically when $n$ is odd.

Then, the main result of \cite{Sa0} for the case
of $n$ even, should be stated as follows.

\begin{thm}
Suppose that $M$ is a closed $n$--dimensional manifold
with $n \geq 2$ even and that $N$ is an oriented surface.
Let $g : M \to N$ be a generic map. Then,
for any nonempty closed oriented $1$--dimensional
submanifold $L$ of $M$ which is
$\Z$--homologous to $S(g)$, there exists
a generic map $f : M \to N$ homotopic to $g$
such that $S(f) = L$.
\end{thm}

\begin{proof}
The proof is almost the same as that given in \cite{Sa0};
however, we need to be careful with the orientations
of $S(g)$ and $L$.

Recall the (oriented) band operation introduced in
\cite[Definition~3.6]{Sa0}.
As the proof of \cite[Lemma~3.9]{Sa0} shows,
$L$ is oriented isotopic to a
$1$--dimensional oriented submanifold of $M$
obtained from $S(g)$ by a finite iteration of 
\emph{oriented} band
operations.

In the proof of \cite[Theorem~2.2]{Sa0} for $n$ even,
just after \cite[Figure~12]{Sa0}, we chose
a joining curve $\lambda$ connecting $p_1$ or $p_2$
and another cusp in $S(g)$.
In fact, both of $p_1$ and $p_2$ can be used, and there
is no difference. This is because we can easily observe that
the oriented band operations associated with $\lambda$ for $p_1$
and for $p_2$ give the same result, considering
the orientations.
Then, in order to adjust the information
mentioned in \cite[Remark~3.11 (2)]{Sa0},
it has been claimed that choosing $p_1$ or $p_2$
appropriately, one can adjust it; however, as
they do not have any difference, this does not work,
unfortunately. Therefore, we need the information
on the orientations. If the band is already oriented,
then the relevant data is automatically adjusted.

Then, the rest of the proof works well without changing
anything. This completes the proof.
\end{proof}

\section*{Acknowledgment}\label{ack}
The author would like to thank Professor \.{I}nan\c{c} Baykur
for stimulating discussions which motivated the
theme of this paper. 
He would also like to thank Dr.\ Rustam Sadykov
for some important comments about certain moves
and for essential discussions.
He would also like to thank Dr.\ Dominik Wrazidlo for some
discussions.
The author would also like to thank Professor Vincent
Blanl\oe il and the Advanced Mathematical Research Institute
of the University of Strasbourg, and also to
the Department of Mathematics, Kansas State University,
for the hospitality
during the preparation of the manuscript.
This work has been supported in part by JSPS KAKENHI Grant Numbers 
JP22K18267, JP23H05437. 
This work was also supported by the Research  
Institute for Mathematical Sciences, an International Joint
Usage/Research Center located in Kyoto University.

\end{document}